\documentclass{cocv}
\usepackage{amsmath,amsfonts,amssymb,amsthm,epsfig}

\usepackage{enumitem}
\usepackage{hyperref}
\usepackage{xcolor}

\numberwithin{equation}{section}

\def\A{{\mathcal A}}

\def\B{{\mathcal B}}
\def\BV{{BV}}
\def\SBVp{{{SBV}^p}}

\def\GSBD{{GSBD}}

\def\GSBDp{{{GSBD}^p}}

\def\div{{\rm div}}
\def\eps{{\varepsilon}}

\def\esp #1 {{{\rm e}^ #1}}
\def\F{{\mathcal F}}

\def\H{{\mathcal H}}

\def\\R^d{{\mathcal \R^d}}
\def\L{{\mathcal L}}

\def\medint{-\kern  -,407cm\int}
\def\medintp{-\kern  -,435cm\int}
\newcommand{\matricim}{\mathbb{M}^{m\times d}}

\def\medintinrigo{-\kern  -,345cm\int}
\def\N{\mathbb{N}}

\def\Om{{\Omega}}

\def\potcap{{\varphi}}

\def\R{\mathbb{R}}
\def\RN{{{\R}^d}}

\def\ristretto{\lfloor}

\def\theta{\vartheta}
\def\d{\mathrm{d}}

\def\tr{\rm tr}

\def\sfera{{\mathchoice
{\setbox0=\hbox{$\displaystyle     \rm S$}\hbox{\raise0.5\ht0\hbox
to0pt{\kern0.35\wd0\vrule height0.45\ht0\hss}\hbox
to0pt{\kern0.55\wd0\vrule height0.5\ht0\hss}\box0}}
{\setbox0=\hbox{$\textstyle        \rm S$}\hbox{\raise0.5\ht0\hbox
to0pt{\kern0.35\wd0\vrule height0.45\ht0\hss}\hbox
to0pt{\kern0.55\wd0\vrule height0.5\ht0\hss}\box0}}
{\setbox0=\hbox{$\scriptstyle      \rm S$}\hbox{\raise0.5\ht0\hbox
to0pt{\kern0.35\wd0\vrule height0.45\ht0\hss}\raise0.05\ht0\hbox
to0pt{\kern0.5\wd0\vrule height0.45\ht0\hss}\box0}}
{\setbox0=\hbox{$\scriptscriptstyle\rm S$}\hbox{\raise0.5\ht0\hbox
to0pt{\kern0.4\wd0\vrule height0.45\ht0\hss}\raise0.05\ht0\hbox
to0pt{\kern0.55\wd0\vrule height0.45\ht0\hss}\box0}}}}
\def\q{{\mathchoice {\setbox0=\hbox{$\displaystyle\rm
Q$}\hbox{\raise
0.15\ht0\hbox to0pt{\kern0.4\wd0\vrule height0.8\ht0\hss}\box0}}
{\setbox0=\hbox{$\textstyle\rm Q$}\hbox{\raise
0.15\ht0\hbox to0pt{\kern0.4\wd0\vrule height0.8\ht0\hss}\box0}}
{\setbox0=\hbox{$\scriptstyle\rm Q$}\hbox{\raise
0.15\ht0\hbox to0pt{\kern0.4\wd0\vrule height0.7\ht0\hss}\box0}}
{\setbox0=\hbox{$\scriptscriptstyle\rm Q$}\hbox{\raise
0.15\ht0\hbox to0pt{\kern0.4\wd0\vrule height0.7\ht0\hss}\box0}}}}

\def\sqr#1#2{\vbox{
         \hrule height .#2pt
         \hbox{\vrule width .#2pt height #1pt \kern #1pt
            \vrule width .#2pt}
         \hrule height .#2pt }}

\newcommand{\res}{\mathop{\hbox{\vrule height 7pt width .5pt depth 0pt
\vrule height .5pt width 6pt depth 0pt}}\nolimits} 

 \providecommand{\Xint}[1]{\mathchoice
    {\XXint\displaystyle\textstyle{#1}}%
    {\XXint\textstyle\scriptstyle{#1}}%
    {\XXint\scriptstyle\scriptscriptstyle{#1}}%
    {\XXint\scriptscriptstyle\scriptscriptstyle{#1}}%
    \!\int}
  \providecommand{\XXint}[3]{{\setbox0=\hbox{$#1{#2#3}{\int}$}
      \vcenter{\hbox{$#2#3$}}\kern-.5\wd0}}
  
  \providecommand{\dashint}{\mathop{\Xint-}}

\begin{document}

\title  [Lower semicontinuity in $GSBD$ for nonautonomous surface integrals]
{Lower semicontinuity in $GSBD$
\\ for nonautonomous surface integrals}

\author{Virginia De Cicco} \address{Dipartimento di Scienze di Base ed Applicate per l'Ingegneria,
Universit\`a di Roma ``La Sapienza",
Via A. Scarpa 16, 00161 Roma, Italy; \email{virginia.decicco@uniroma1.it }}
\author{Giovanni Scilla} \address{Dipartimento di Matematica ed Applicazioni "R. Caccioppoli", 
   Universit\`{a} degli Studi di Napoli Federico II, 
Via Cintia, Monte Sant'Angelo 80126 Napoli, Italy; \email{giovanni.scilla@unina.it} } 

\begin{abstract} 
We provide a sufficient condition for lower semicontinuity of nonautonomous noncoercive surface energies defined on the space of $GSBD^p$ functions, 
whose dependence on the $x$-variable is $W^{1,1}$ or even $BV$: the notion of \emph{nonautonomous symmetric joint convexity}, which extends the analogous definition devised for autonomous integrands in \cite{FPS}  where the conservativeness of the approximating vector fields is assumed. This condition allows to extend to our setting a nonautonomous chain formula in $SBV$ obtained in \cite{ACDD}, and this is a key tool in the proof of the lower semicontinuity result. This new joint convexity can be checked explicitly for some classes of surface energies arising from variational models of fractures in inhomogeneous materials. 
\end{abstract}



\maketitle

\section{Introduction}\label{s:introduction}

The modern variational approach to the quasistatic growth of brittle cracks in elastic materials, proposed by Francfort and Marigo~\cite{FM} 
following on Griffith's theory~\cite{Gr}, assumes that at each time the equilibrium of a crack is the result of the balance between the bulk elastic energy released (due to the crack growth) and the surface energy dissipated to produce a new portion of crack. 
In the formulation of the related variational problems, the energy dissipated to produce the crack 
can be described by a surface integral of the type
\begin{equation}
{\mathcal  K}(u):=\int_{\Omega\cap J_u} \phi(x,u^-,u^+,\nu_{u})\,\mathrm{d}\H^{d-1},
\label{eq:surfaceterm}
\end{equation}
where $\Omega$ denotes the reference configuration, $u$ the displacement field, $J_u$ the crack surface and $\H^{d-1}$ denotes the $(d-1)$-dimensional Hausdorff measure. In the case of a constant density in \eqref{eq:surfaceterm}, $\phi$ coincides with the ``toughness'' of the material and ${\mathcal  K}(u)$ is proportional to the measure of the crack $\mathcal{H}^{d-1}(J_u)$. More in general, $\phi$ may depend on the crack opening $[u]:=u^+-u^-$, where $u^+$ and $u^-$ are the traces of $u$ on both sides of its jump set $J_u$, thus modeling a cohesive or Barenblatt type fracture (see~\cite{B}). The anisotropies and inhomogeneities in the body are taken into account through the possible dependence of $\phi$ on the normal $\nu_u$ to $J_u$ and the position $x$. 
The bulk elastic energy released during the crack formation is modeled by a functional of the form
\begin{equation}
\mathcal{W}(u):=\int_\Om W(x,\nabla u)\,\mathrm{d}x\,,
\label{eq:bulken}
\end{equation}
where $\nabla u$ is the deformation gradient and the potential $W(x,\xi)$ is typically assumed to be a Carath\'eodory function; i.e., measurable in $x$ and continuous in $\xi$ for a.e. $x$, such that $W(x,\cdot)$ is convex and $W(x,\xi)\geq c|\xi|^p$, $p>1$, for a.e. $x\in\Omega$.

A key ingredient to prove the existence of minimizers for functionals of the type
\begin{equation} 
u\to\mathcal{E}(u):=\mathcal{W}(u) + {\mathcal  K}(u) 
\label{eq:minprob}
\end{equation}
is the lower semicontinuity of surface integrals \eqref{eq:surfaceterm}. Within the weak theory in the space $(G)SBV$ of (generalized) special functions of bounded variation, the problem of finding sufficient conditions for the lower semicontinuity of functionals like \eqref{eq:surfaceterm} has been addressed e.g. in \cite{ADCF3,Amb2, AmbBra1,AFP, DMFT,DT, DC}. One of these conditions is the \emph{$BV$-ellipticity} of the integrand, introduced by Ambrosio and Braides in \cite{AmbBra1}. More recently, aiming to describe the evolution in plasticity and elastoplasticity and the formation of microstructures, the space $\GSBD$ of ``generalized special functions of bounded deformation", introduced in \cite{DMJems}, has been considered as the proper functional space. Many progresses in this direction, in particular concerning the study of free-discontinuity problems, can be found in the recent literature, see \cite{CCI, CC, CC2, CFI2, Cri, CF, FPS, FS1, FS2}.  

The scope of this paper is to study the lower semicontinuity in $GSBD$ of nonautonomous surface integrals of the type \eqref{eq:surfaceterm} where the surface integrand $\phi$ depends explicitly and possibly in a discontinuous way on the space variable $x$, and without any coercivity assumption. Our contribution is then motivated by obtaining the well posedness on $GSBD^p$ (here $p$ stands for the exponent of summability of $e(u):=\frac{1}{2}(\nabla u + \nabla u^T)$, the approximate symmetrized gradient) of variational problems associated to energies $\mathcal{E}(u)$. 

The existence of minimizers for the functional \eqref{eq:minprob} is guaranteed in $SBD^p$ by the
compactness result \cite[Theorem~1.1]{BCDM}, provided one has an a priori bound for $u$
in $L^\infty$. Unfortunately, it is hard to obtain such a bound, which is in fact quite unrealistic in Fracture Mechanics. 
To overcome this difficulty, the landing in the space $GSBD^p$ of the generalised $SBD^p$ functions is needed. Our main goal is to prove a lower semicontinuity result in $GSBD^p$ for functionals as \eqref{eq:minprob}, along sequences $\{u_n\}$ in $\GSBDp(\Om;\R^d)$ such that $u_n\to u$ in measure on $\Om$ and $\Vert e (u_n)\Vert_p$, $\H^{d-1}(J_{u_n})$ are uniformly bounded with respect to $n\in\N$. Indeed, the bulk energy \eqref{eq:bulken} here depends only on the symmetric part of $\nabla u$. 
It is worth mentioning that we have to consider in \eqref{eq:minprob} the additional term $\int_\Om \Psi(|u|)\,\mathrm{d}x$, for a continuous $\Psi$ such that $\displaystyle\lim_{s\to+\infty}\Psi(s)=+\infty$. This is only needed to apply the compactness result in $GSBD^p$ (see \cite[Theorem~11.3]{DMJems}) and, actually, should be dropped in favor of Dirichlet boundary conditions, provided we agree to work in the larger space $GSBD^p_\infty$, recently introduced in \cite{CF}. Since the lower semicontinuity of this additional term is trivial, and the bulk term $\mathcal{W}(u)$ is lower semicontinuous by convexity, 
only the problem of the lower semicontinuity of the surface term \eqref{eq:surfaceterm} has to be addressed. 

In the autonomous case a necessary and sufficient condition for the lower semicontinuity is the \emph{$BD$-ellipticity} introduced in \cite{FPS}, by adapting to the $BD$ setting the notion of $BV$-ellipticity. 
Both the conditions are very difficult to check directly.
If on the one hand, a sufficient condition for the $BV$-ellipticity is the \emph{joint convexity}; i.e., the existence of continuous vector fields $g_j:\R^d\to \R^d$ such that 
$$
\phi(r,t,\xi)=\sup_{j\in\N} \langle g_j(r)-g_j(t),\xi\rangle\quad\quad
{\rm for\ all\ } (r,t,\xi)\in  \R^d\times \R^d\times  \R^{d}\,,
$$
on the other in \cite{FPS} the notion of \emph{symmetric} joint convexity was introduced as a sufficient condition for $BD$-ellipticity, by requiring the further condition that  for every $j\in\N$  the vector fields $g_j$ are \emph{conservative}. Roughly speaking, this latter ensures that an integration by part formula holds in the $BD$ context.  Nevertheless, while the joint convexity (and the $BV$-ellipticity) 
are nowadays well-understood for both autonomous and nonautonomous surface integrals (see \cite{AmbBra1, AFP,DC}), the analogous problems in the modern setting of $GSBD$ are still largely open. The few examples of lower semicontinuous functionals in $BD$ present in literature (see \cite{CC22, DOT, FPS, GZ1, KP} and the references therein) are mostly concerning with the autonomous setting. In the rare cases where an explicit dependence on the spatial variable is allowed, the integrands are assumed to be continuous in $x$ (see, for instance, \cite{CC2}), or a coercivity assumption in the gradient variable is required, as in \cite{KP} 
where $\phi(x,\cdot)$ is a Finsler metric. 

\emph{Our results:} 
we address the problem of the lower semicontinuity for surface integrands $\phi$ depending explicitly on the spatial variable $x$. The new feature is that $\phi$ may possibly be discontinuous with respect to $x$, but in a ``controlled'' way, as it admits a $BV$-regularity. 
We also renounce to any coercivity assumption, so that finding sufficient conditions for the lower semicontinuity is non trivial and involves the regularity of $\phi$ at $x$. Firstly, we assume a $W^{1,1}$-dependence in $x$ and introduce the notion of {\it nonautonomous ($NA$) symmetric jointly convex} integrand. Namely, we require the existence of a sequence of nonautonomous vector fields $g_j:\Omega\times  \R^d\to \R^d$ such that 
\begin{equation}
\phi(x,r,t,\xi)=\sup_{j\in\N} \langle g_j(x,r)-g_j(x,t),\xi\rangle\quad\quad
{\rm for\ all\ } (x,r,t,\xi)\in \Omega\times \R^d\times \R^d\times  \R^{d}
\label{eq:nasymmjointconv}
\end{equation}
and each $g_j$ satisfies some conditions including the $W^{1,1}$ dependence w.r.t. $x$ and the conservativeness, see \ref{ass-a} and \ref{ass-f}.

In order to prove the lower semicontinuity result for the corresponding functionals, a crucial tool is a nonautonomous chain rule formula (proven in Theorem~\ref{chain rule1}) for the divergence of the composition $v(x)=g_j(x,u(x))$ with $u\in GSBD^p$.
This formula is obtained by combining the analogous in $SBV$ \cite{ACDD} with a recent approximation result due to Chambolle and Crismale (see \cite{CC}), and extends to the nonautonomous setting the integration by parts formula provided by \cite[Lemma~5.3]{FPS}, and to the $GSBD^p$ setting the nonautonomous chain rule \cite[Theorem~2.2]{ACDD} itself. Section~\ref{sec: joint} is devoted to the proof of main result, the lower semicontinuity theorem (Theorem~\ref{thm:gsbdsemicontinuity1}): every surface integral with a nonautonomous symmetric jointly convex integrand is lower semicontinuous on $GSBD^p$. We first prove the assertion for those integrands that can be represented as in \eqref{eq:nasymmjointconv} with vector fields $g_j(x,r)$ satisfying all the assumptions \ref{ass-a}--\ref{ass-f} of the nonautonomous chain rule, Theorem~\ref{chain rule1}. Then we obtain the general case in a standard way via regularization of the approximating vector fields. 

As a next step, in Section~\ref{sec: bvjoint}, we relax all our results to allow the functions for a more general $BV$-dependence on $x$, by introducing the notion of {\it $ {{BV}}$ symmetric jointly convex} function. Note that this is a generalization of the $NA$ symmetric joint convexity only for certain integrands (cf. the discussion at the beginning of Section~\ref{sec: bvjoint}). In particular, as a nontrivial example, in Proposition~\ref{thm:thm4.13fpsbv} we prove the $BV$ symmetric joint convexity for integrands of the form
\begin{equation}
\phi(x,r,t,\xi) := \kappa(x,\xi)\,,
\label{eq:separatedint}
\end{equation}
where $\kappa(\cdot,\xi)\in BV$ (cf. \ref{ass-k1cc}) and $\kappa(x,\cdot)$ is even, positively 1-homogeneous, and convex for a.e. $x$. Despite of its simple structure, this example already contains some features related to the nonautonomous setting. Indeed, the explicit construction of the vector fields $g_j(x,r)$ relies on a classical approximation result for real-valued convex functions due to De Giorgi (see Lemma~\ref{degiorgi}), which provides an explicit formula for the approximating affine functions.

We rediscover some approximation arguments developed in \cite{DC} to obtain a lower semicontinuity result in $GSBD^p$ for  $ {{BV}}$ symmetric jointly convex integrands, see Proposition~\ref{t:SCI1}. The key point there is the proof of the result for ``splitting-type'' integrands of the form
\begin{equation}
\phi(x,r,t,\xi) :=a(x)\kappa(r,t,\xi)\,,
\label{eq:splitfunctintro}
\end{equation}
where $a$ is a $BV$ function and $\kappa$ is symmetric jointly convex (Proposition~\ref{corcor}). Indeed, $a$ can be approximated from below by $W^{1,1}$ functions, and Theorem~\ref{thm:gsbdsemicontinuity1} can be applied to the approximating $NA$ symmetric jointly convex functions. Then, again an approximation argument by means of functions as in \eqref{eq:splitfunctintro} gives the lower semicontinuity for functionals whose integrands are strictly positive and $BV$ symmetric jointly convex. Eventually, in Section~\ref{sec: existence} we apply our lower semicontinuity results to prove the existence of minimizers for functionals of the type \eqref{eq:minprob}.

We conclude the presentation of our results mentioning that a further issue to be addressed is the introduction of a notion of ``nonautonomous $BD$-ellipticity'' together with the investigation of its connections with \emph{NA}/$BV$ symmetric joint convexity. This problem seems to be as interesting as challenging with the available mathematical tools, and therefore has to be deferred to future studies.

\section{Preliminaries}\label{sec:preliminaries}

\subsection{Notation}
Let $\Om$ be a bounded open subset of $\RN$. 
We denote by $\A(\Om)$ the family of all open subsets $A$
of $\Om$ and by $\B(\Om)$ the $\sigma$-algebra of all Borel subsets
$B$ of $\Om$. For a set $E\subset\Omega$, we will often denote by $E^c$ its complement $\Omega\backslash E$.

\noindent
Let ${\L}^{d}$ denote the Lebesgue measure on $\RN$
and ${\H}^{d-1}$ the Hausdorff
measure of dimension $(d-1)$ on $\RN$. We denote by $\langle \cdot,\cdot\rangle$ the standard Euclidean scalar product on $\R^d$. $\mathbb{S}^{d-1}$ is the unit sphere. We denote by $\mathbb{M}^{m\times n}$ the set of real $m\times n$ matrices, and by $\mathbb{M}_{\rm sym}^{n\times n}$ the set of all real symmetric $n\times n$ matrices. Given a matrix $A$, $A^T$ stands for the  transpose of $A$ and ${\rm tr}(A)$ for its trace. Given two matrices $A,B\in\mathbb{M}^{m\times n}$, the Frobenius scalar product will be denoted by $A:B:={\rm tr}(A^TB)$, while $|A|:=\sqrt{{\rm tr}(A^TA)}$ will indicate the associated norm. For $x\in\Omega$ and $\rho>0$, we denote by $B_\rho(x)$ the open ball centred at $x$ with radius $\rho$. We will denote by $\mathcal{C}(X;Y)$ the space of continuous functions from $X$ to $Y$, while by $C^1_c(X)$ and $C^\infty_c(X)$ the spaces of $C^1$ and $C^\infty$ functions with compact support on $X$, respectively. The symbols $L^p(X)$ and $W^{k,p}(X)$ stand for the classical Lebesgue and Sobolev spaces defined on $X$, respectively.
\medskip

\subsection{$BV$ and $SBV$ functions.}\label{s:bv}

For a general survey on the spaces of $\BV$ and $SBV$ functions 
we refer for instance to \cite{AFP}. Below, we just recall some basic definitions useful in the sequel. 

If  $u\in L^1_{\rm loc}(\Om;\R^m)$ and $x\in\Omega$, the {\it precise representative of $u$ at $x$} is defined
as the unique value $\widetilde{u}(x)\in\R^m$ such that
$$
\lim_{\rho\to 0^+} \frac{1}{\rho^d}\int_{{B_\rho(x)}}\!|u(y)-\widetilde{u}(x)|\,\d x=0\,.
$$
The set of points in $\Omega$ where the precise representative of $x$ is not defined is called the
{\it approximate singular set} of $u$ and denoted by $S_u$. We say that a point $x\in\Omega$ is an approximate jump point of $u$
if there exist $a,b\in\R^m$ and $\nu\in\mathbb{S}^{d-1}$, such that $a\not = b$ and
$$
\lim_{\rho\to 0^+}\dashint_{B^+_\rho(x,\nu)} |u(y)-a|\, \d y=0
\qquad{\rm and}\qquad
\lim_{\rho\to 0^+}\dashint_{B^-_\rho(x,\nu)} |u(y)-b|\, \d y=0
$$
where $B^\pm_\rho(x,\nu):= \{y\in B_\rho(x)\ :\ \langle y-x,\nu\rangle\gtrless0\}$.
The triplet $(a,b,\nu)$ is uniquely determined by the previous formulas, up to a permutation
of $a,b$ and a change of sign of $\nu$, and it is denoted by $(u^+(x),u^-(x),\nu_u(x))$.
The Borel functions $u^+$ and $u^-$  are called the {\it upper and
lower approximate limit} of $u$ at the point $x\in\Omega$. The set of approximate jump points of $u$ is denoted by $J_u$.

The space ${BV}(\Om;\R^m)$ of {\it functions of bounded variation} is defined as the set
of all $u\in L^1(\Om;\R^m)$ whose distributional gradient $Du$ is a bounded Radon measure on $\Om$
with values in the space $\matricim$ of $m\times d$ matrices. Moreover, the usual decomposition

\begin{equation*}
Du = \nabla u\,\L^d + D^c u + (u^+-u^-)\otimes \nu_u\,\H^{d-1}\ristretto{J_u}
\end{equation*}
\noindent
holds, where $\nabla u$ is the Radon-Nikod\'ym derivative of $Du$ with respect to the Lebesgue measure and $D^cu$ is  the {\it Cantor part} of $Du$. For the sake of simplicity, we denote by $D^su =D^c u + (u^+-u^-)\otimes \nu_u\,\H^{d-1}\ristretto{J_u}$.

We recall that the space ${SBV}(\Om;\R^m)$ of {\it special functions of bounded variation} is defined as the set
of all $u\in \BV(\Om;\R^m)$ such that $D^su$ is concentrated on $S_u$; i.e., $|D^su|(\Om\setminus S_u)=0$. Finally, for $p>1$ the space $\SBVp(\Om;\R^m)$ is the set of $u\in SBV(\Om;\R^m)$ with $\nabla u\in L^p(\Om;\matricim)$ and $\H^{d-1}(S_u)<\infty$.

\subsection{$GBD$ and $GSBD$ functions}\label{sec: gsbd}

In this section we recall some basic definitions and results on generalized functions with bounded deformation, as introduced in \cite{DMJems}. Throughout the paper we will use standard notations for the space $(G)SBD$, referring the reader to \cite{ACDM, BCDM, Temam} for a detailed treatment on the topic.

Let $\xi\in\R^d\backslash\{0\}$ and $\Pi^\xi=\{y\in\R^d:\, \langle\xi,y\rangle=0\}$. If $\Omega\subset\R^d$ and $y\in\Pi^\xi$ we set $\Omega_{\xi,y}:=\{t\in\R:\, y+t\xi\in \Omega\}$. 
Given $u:\Omega\to\R^d$, $d\geq2$, we define $u^{\xi,y}: \Omega_{\xi,y}\to\R$ by $u^{\xi,y}(t):=\langle u(y+t\xi),\xi\rangle$. 

\noindent
We then have the following definitions:

$(i)$ An $\mathcal L^{d}$-measurable function $u:\Omega\to \R^{d}$ belongs to $GBD(\Omega;\R^d)$ if there exists a positive bounded Radon measure $\lambda_u$ such that, for all $\tau \in C^{1}(\R^{d})$ with $-\frac12 \le \tau \le \frac12$ and $0\le \tau'\le 1$, and all $\xi \in \mathbb{S}^{d-1}$, the distributional derivative $D_\xi (\tau(\langle u,\xi\rangle))$ is a bounded Radon measure on $\Omega$ whose total variation satisfies
$$
\left|D_\xi (\tau(\langle u,\xi\rangle))\right|(B)\le \lambda_u(B)
$$
for every Borel subset $B$ of $\Omega$. 

$(ii)$ A function $u \in GBD(\Omega;\R^d)$ belongs to the subset $GSBD(\Omega;\R^d)$ of special functions of bounded deformation if, in addition, for every $\xi \in \mathbb{S}^{d-1}$ and $\mathcal H^{d-1}$-a.e.\ $y \in \Pi^\xi$, it holds that $u^{\xi,y}\in SBV_{\mathrm{loc}}(\Omega_{\xi,y})$.

We recall that the inclusions $BD(\Omega;\R^d)\subset GBD(\Omega;\R^d)$ and $SBD(\Omega;\R^d)\subset GSBD(\Omega;\R^d)$ hold (see \cite[Remark 4.5]{DMJems}).  {Although} they are, in general, strict,  relevant properties of $BD$ functions are retained also in this weak setting. In particular, $GBD$-functions have an approximate symmetric differential $e(u)(x)$ at $\mathcal L^{d}$-a.e.\ $x\in \Omega$. Furthermore the jump set $J_u$ of a $GBD$-function is $\mathcal H^{d-1}$-rectifiable (this is proven in \cite[Theorem 6.2 and Theorem 9.1]{DMJems}).

{Let $p>1$.} The space $GSBD^p(\Omega;\R^d)$ is defined as
$$
GSBD^p (\Omega;\R^d):= \{u \in GSBD(\Omega;\R^d):\, e(u) \in L^p (\Omega; \mathbb M_{\mathrm{sym}}^{d\times d})\,,\,\mathcal H^{d-1}(J_u) < +\infty\}\,.
$$

{ Given a functional $\mathcal{E}: GSBD^p(\Omega;\R^d)\to[0,+\infty]$, we say that $\mathcal{E}$ is \emph{lower semicontinuous in $GSBD^p$} if for every $(u_n)_n\subset GSBD^p(\Om;\R^d)$ converging in measure to $u\in GSBD^p(\Om;\R^d)$ 
and such that
\begin{equation*}
\sup_{n\in\N}\left[\int_\Om |e(u_n)|^p\,\mathrm{d}x + \H^{d-1}(J_{u_n})\right] <+\infty\,,
\end{equation*}
we have
\begin{equation*}
\mathcal{E}(u) \leq
\liminf_{n\to +\infty}\mathcal{E}(u_n)\,.
\end{equation*}
}

\subsection{Approximation results}\label{approx}

Let us recall some well known approximation results.

The first one is very general and concerns the lower semicontinuity of a functional whose integrand is the supremum of a sequence of 
integrands corresponding to lower semicontinuous functionals. For our purposes, we formulate the statement in $GSBD^p$.  

\begin{lmm}\label{l:passup}
Let $h,h_j:\Om\times\R^d\times\R^d\times\RN \to [0,+\infty)$,  $j\in\N$, be Borel functions such that 
\begin{itemize}
\item[$(i)$] the functionals $\F_{h_j}$ defined by
$$
\F_{h_j}(u,\Omega):= \int_{\Omega\cap J_u} h_j(x,u^-,u^+,\nu_{u})\,\d\H^{d-1}
$$
are 
lower semicontinuous in $GSBD^p(\Om;\R^d)$ for every $j\in\N$; 
\item[$(ii)$] $h(x,r,t,\xi) = \sup_{j\in\N} h_j(x,r,t,\xi)$ for all $(x,r,t,\xi)\in(\Om\setminus N)\times\R^d\times\R^d\times\RN,$
where $N\subset\Om$ is a Borel set with $\H^{d-1}(N)=0$.
\end{itemize}
Then, also
$$
\F_{h}(u,\Omega):= \int_{\Omega\cap J_u} h(x,u^-,u^+,\nu_{u})\,\d\H^{d-1}
$$
is  lower semicontinuous in $GSBD^p(\Om;\R^d)$.
\end{lmm}
\begin{proof}
The argument is nowadays standard, but for the reader's convenience we prefer to give a proof (see, e.g., Step~2 of \cite[Proof of Corollary~2.6]{FPS} for a similar result). Let $(u_n)_{n\in\N}\subset GSBD^p(\Omega;\R^d)$ and $u\in GSBD^p(\Omega;\R^d)$ be such that $u_n\to u$ in measure on $\Omega$ as $n\to+\infty$. Let $V\subset\Omega$ be any open set, and define the non-negative superadditive function
\begin{equation*}
\Lambda(V):=\mathop{\lim\inf}_{n\in\N} \F_{h}(u_n,V)\,.
\end{equation*}
Since each $\F_{h_j}$ is lower semicontinuous, we have $\Lambda(V)\geq \F_{h_j}(u,V)$ for every $j\in\N$ and every open set $V\subset\Omega$. Then, by a classical lemma on the supremum of measures (see, e.g., \cite[Lemma~2.8]{FPS}), we get
\begin{equation*}
\Lambda(V)\geq  \int_{V\cap J_u} \sup_{j\in\N} h_j(x,u^-,u^+,\nu_{u})\,\d\H^{d-1} = \F_{h}(u,V),
\end{equation*}
whence the desired assertion follows choosing $V=\Omega$.
\end{proof}

The second one is an approximation result for $GSBD^p$ functions, stated in \cite[Theorem~1.1]{CC}, which we recall here in a slightly simplified version.

\begin{thrm}\label{thm:cc}
Let $u\in GSBD^p(\Omega;\R^d)$, $p > 1$. Then, there exists a sequence of functions $(u_k)_k\subset SBV^p(\Omega;\R^d)\cap L^\infty(\Omega;\R^d)$ such that each $J_{u_k}$ is closed in $\Omega$ and included in a finite union of closed connected pieces of $C^1$ hypersurfaces, $u_k\in W^{1,\infty}(\Omega\backslash J_{u_k}; \R^d)$, and
\begin{enumerate}
\item[(i)] $u_k\to u$ a.e. on $\Omega$;
\item[(ii)] $\|e(u_k)-e(u)\|_{L^p(\Omega)}\to0$;
\item[(iii)] $\mathcal{H}^{d-1}(J_{u_k}\triangle J_u)\to0$;
\item[(iv)] $\int_{J_{u_k}\cup J_u}\tau(|u_k^\pm-u^\pm|)\,\mathrm{d}\mathcal{H}^{d-1}\to0$,
\end{enumerate}
for some $\tau\in C^1(\R)$ with $-\frac{1}{2}\leq \tau \leq \frac{1}{2}$, $0\leq \tau' \leq1$, and $\{\tau=0\}=\{0\}$.
\end{thrm}

We recall now a compactness result in $GSBD^p$ (see \cite[Theorem~11.3]{DMJems}).

\begin{thrm}\label{thm:compactness}
Let $(u_k)_k$ be a sequence in $GSBD^p(\Omega;\R^d)$. Assume that there exist a constant $M>0$ and an increasing continuous function $\Psi:\R^+\to\R^+$ with $\displaystyle\lim_{s\to+\infty}\Psi(s)=+\infty$ such that
\begin{equation}
\int_\Om \Psi(|u_k|)\,\mathrm{d}x + \int_\Omega |e(u_k)|^p\,\mathrm{d}x + \mathcal{H}^{d-1}(J_{u_k})\leq M
\label{eq:equicoerciv}
\end{equation}
for every $k$. Then there exist a subsequence (not relabeled) and a function $u\in GSBD^p(\Omega;\R^d)$ such that
\begin{equation*}
\begin{split}
& u_k\to u \quad \mbox{ pointwise $\mathcal{L}^d$-a.e. on $\Omega$,}\\
& e(u_k)\rightharpoonup e(u)  \quad \mbox{ weakly in $L^p(\Omega;\R_{\rm sym}^{d\times d})$,} \\
&\mathcal{H}^{d-1}(J_u) \leq \mathop{\lim\inf}_{k\to+\infty} \mathcal{H}^{d-1}(J_{u_k})\,.
\end{split}
\end{equation*}
\end{thrm}

We conclude this section with a classical approximation result {for} convex functions by means of affine functions due to De
Giorgi~\cite{deg} (see also \cite[Theorem~4.79]{FL}). Let $f:\Om\times \R^d\to [0,+\infty)$ be convex and positively $1$-homogeneous in the last variable, and let $\alpha\in C^1_c(\R^d)$ be a non negative function such that $\int_{\R^{d}}
\alpha(\xi)\,\mathrm{d}\xi=1$. 
For every $j\in\N$ and $q\in \mathbb{Q}^d$, we define
\begin{align}
a^0_{j,q}(x) &  =\int_{\R^{d}}f(x,\xi)\Big(  (d+1)\alpha_{j,q}
(\xi)+\langle\nabla\alpha_{j,q}(\xi),\xi\rangle\Big)  \,\mathrm{d}\xi\label{coeff}\\
a_{j,q} (x) & =-\int_{\R^{d}}f(x,\xi)\nabla\alpha_{j,q}(\xi)\,\mathrm{d}\xi\,,\label{coeff1}
\end{align}
where $\alpha_{j,q}(\xi):=j^d\alpha(j(q-\xi))$. We then have the following approximation of $f$ from below with a sequence of affine functions.

\begin{lmm}\label{degiorgi}
Let $f:\Om\times \R^d\to [0,+\infty)$ be convex 
in the last variable, and define for every $x\in\Omega$ the sequence $(a_{j,q}(x))_{j,q}$, $j\in\N$, $q\in\mathbb{Q}^d$ as in \eqref{coeff} and \eqref{coeff1}. Then for all $(x,\xi)\in\Omega\times\R^d$ we have
$$
f(x,\xi)=\sup_{j\in\N\,,\,\,q\in\mathbb{Q}^d}\, [a^0_{j,q}(x)+\langle a_{j,q}(x),\xi\rangle]^+\,.
$$
If $f(x,\cdot)$ is also positively $1$-homogeneous,
\begin{equation}\label{poppop1}
f(x,\xi)=\sup_{j\in\N\,,\,\,q\in\mathbb{Q}^d}\langle a_{j,q}(x),\xi\rangle^+\,.
\end{equation}
\end{lmm}

The approximation \eqref{poppop1} is very useful in semicontinuity problems since the coefficients $a_{j,q}(x)$ above depend explicitly on $f$ and on the lower order variable $x$. Thus, the regularity properties of $f(\cdot,\xi)$ are inherithed by $a_{j,q}(\cdot)$ through formulas \eqref{coeff1} and \eqref{poppop1}. Note also that the explicit dependence of $a_{j,q}$ on $q$ could be neglected by considering in \eqref{coeff1} a vector $q_j\in\mathbb{Q}^d$ for each $j\in\N$. In this case, we may set $a_j:=a_{j,q_j}$ and compute the supremum in \eqref{poppop1} over $j\in\N$.

\begin{rmrk}\label{rem:widerclass}
Notice that if $f(x,\cdot)$ is even, convex, positively $1$-homogeneous and bounded away from zero on $\mathbb{S}^{d-1}$, then
for every fixed $x$ it is the support function of the symmetric compact convex set with nonempty interior $K(x)$ given by
\begin{equation*}
K(x):=\left\{z\in\R^d:\,\, \langle z,\xi\rangle \leq f(x,\xi)\,\,\mbox{ for every $\xi\in\R^d$}\right\}\,.
\end{equation*}
Then, by \eqref{poppop1}, for every fixed $x\in\Omega$, we have $a_j(x)\in K(x)$ for every $j\in\N$. Equivalently, the function $a_j:\Omega\to\R^d$ is a selection of the multifunction $K:\Omega\to \mathcal{P}(\R^d)$ for every $j\in\N$ {(see \cite[Chapter~III]{CastVal})}, where $\mathcal{P}(\R^d)$ denotes the power set of $\R^d$. Moreover, if $a_j(x)$ is an interior point of $K(x)$, for every $v\in \mathbb{S}^{d-1}$ we can find $\sigma>0$ small enough such that $a_j(x)+\sigma v\in K(x)$.
For our purposes, it is useful to introduce the (nonempty) countable sets
\begin{equation}
\begin{split}
\mathcal{A}(x):=\left\{a_j(x)+\sigma v:\,\, \sigma\in\mathbb{Q}^+,\,v\in\mathbb{Q}^d\cap \mathbb{S}^{d-1}\,,j\in\N \right\}\cap K(x)\,,\quad x\in\Omega\,. 
\end{split}
\label{eq:classdomain}
\end{equation}
We will denote by $b_l:\Omega\to\R^d$, $l\in\N$ the {measurable} selections of the multifunction $\mathcal{A}:\Omega\to \mathcal{P}(\R^d)$.


\end{rmrk}

\subsection{Capacity}\label{not5}
Following \cite[Section~2.3]{DC}, we briefly recall the notion of 1-capacity of a set and its connections with $BV$ and Sobolev functions.
Given an open set $A\subset\R^d$, the {\it $1$-capacity} of $A$ is defined by setting  \begin{equation*}
C_1(A) : = \inf\left\{\int_{\R^d } |D\potcap|\,dx\ :\ \potcap\in W^{1,1}(\R^d),
\quad \potcap\geq 1\quad \L^{d}{\rm -a.e.\ on}\ A\right\}\,.
\end{equation*}
Then, the $1$-capacity of an arbitrary set $B\subset\R^d$ is given by
\begin{equation*}
C_1(B) := \inf\{C_1(A)\ :\ A\supseteq B,\ A\ {\rm open}\}\,.
\end{equation*}
It is well known 
that for every Borel set $B\subset\R^d$
$$
C_1(B)=0\qquad\Longleftrightarrow\qquad
\H^{d-1}(B)=0\,.
$$


We recall that a function $g:\R^d\to\R$ is said $C_1$-quasi continuous
if for every $\eps>0$ there exists an open set $A$, with $C_1(A)<\eps$, such that
$g|^{}_{A^c}$ is continuous on $A^c$;
$C_1$-quasi lower semicontinuous  and $C_1$-quasi upper semicontinuous
functions are defined similarly.

It is well known that if $g$ is a $W^{1,1}$ function, then its precise representative $\widetilde{g}$ is
$C_1$-quasi continuous (see \cite[Sections 9 and 10]{FZ}).
Moreover, to every $\BV$ function $g$, it is possible to associate a $C_1$-quasi lower semicontinuous
and a $C_1$-quasi upper semicontinuous representative, as stated by the
following theorem (see \cite{CDLP}, Theorem 2.5).

\begin{thrm}\label{t:lecce}
For every function $g\in\BV(\Om)$, the approximate upper limit $g^+$ and the approximate
lower limit $g^-$ are $C_1$-quasi upper semicontinuous and $C_1$-quasi lower semicontinuous,
respectively.
\end{thrm}

Moreover we recall the following lemma which is an approximation result due to Dal Maso (see \cite{DM83},
Lemma 1.5 and \S 6).

\begin{lmm}\label{maso1}
Let $g:\R^d\to[0,+\infty)$ be a $C_1$-quasi  lower semicontinuous function.
Then there exists an increasing sequence of nonnegative
functions $\{g_h\}\subseteq W^{1,1}(\R^d)$ such that, for every $h\in\N$,
$g_h$ is approximately continuous $\H^{d-1}$-almost everywhere in $\R^d$ and
$g_h(x)\to g(x)$ as $h\to+\infty$  for $\H^{d-1}$-almost every $x\in\R^d$.
\end{lmm}


\section{Nonautonomous chain rule formula for the divergence}\label{chaingsbd}

The aim of this section is the proof of a nonautonomous chain rule formula for the divergence of the composition 
\begin{equation*}
v(x):=g(x,u(x))\,,
\end{equation*}
where $u\in GSBD^p$ and
$g\colon\R^d\times\R^d\to\R^d$ complies with 
\leavevmode
\begin{enumerate}[font={\normalfont},label={(G\arabic*)}]
\item $x\mapsto g(x,r)$ belongs to $W^{1,1}_{\rm loc}(\R^d;\R^d)$ for all
$r\in\R^d$, {and there exists a positive function $h_1\in L^1_{\rm{loc}}(\R^d)$ such that $
|g(x,r)|\leq h_1(x)$ for all $r\in\R^d$ and  for $\L^{d}$-a.e. $x\in\R^d$}\,; \label{ass-a}
\item there exist a positive function 
$h_2\in L^1_{\rm loc}(\R^d)$ and a modulus of continuity $\omega:[0,\infty)\to[0,1]$ such that
\[
|\nabla_x g(x,r)-\nabla_x g(x,s)|\le
\omega(|r-s|)h_2(x)
\]
for all $r,\,s\in \R^d$ and  for $\L^{d}$-a.e. $x\in\R^d$\,;\label{ass-b}
\item
there exists a Lebesgue negligible set $N\subset \R^d$
such that $r\mapsto g(x,r)$ is continuously differentiable in $
\R^d$ for all $x\in\R^d\setminus N$\,; \label{ass-c}
\item there exists a constant $M>0$ such that $|\nabla_r g(x,r)|\le M$ for all
$x\in\R^d\setminus N$ and $r\in\R^d$\,; \label{ass-d}
\item 
there exists a modulus of continuity
$\tilde\omega$ independent of $x$ such that
\[
 |\nabla_r g(x,r)-\nabla_r g(x,s)|\leq \tilde \omega(|r-s|)
 \]
for all $r,\,s\in \R^d$ and $x\in\R^d\setminus N$\,. \label{ass-e}

\end{enumerate}

First, we notice that assumptions {\rm\ref{ass-a}}-{\rm\ref{ass-e}} are enough to obtain a vectorial chain rule formula in $SBV$ for nonautonomous functionals, which will be a technical tool for the proof of an analogous result in $GSBD$. We recall it for reader's convenience with the statement below, Theorem~\ref{chain rule}, and it can be seen as a 
simple case of a chain rule in $BV$ proven in \cite[Theorem~2.2]{ACDD} under more general assumptions on the dependence in $x$\,. 

\begin{thrm}\label{chain rule}
Assume that $g$ complies with {\rm\ref{ass-a}}-{\rm\ref{ass-e}} above.
Then there exists a set $\mathcal N\subset\R^d$ with $\H^{d-1}(\mathcal N)=0$\,, such that,
 for every $r\in \R^d$ the function $g(\cdot,r)$ is approximately continuous in $\R^d\setminus \mathcal N$, 
and
for any function $u\in SBV_{\rm loc}(\R^d;\R^d)\cap L^\infty_{\rm loc}(\R^d;\R^d)$,
the function $v(x):=g(x,u(x))$ belongs to $
SBV_{\rm loc}(\R^d;\R^d)$ and the following chain rule holds:
\begin{itemize}
\item[(i)] (Lebesgue part) 
for $\L^{d}$-a.e. $x$ the map $y\mapsto g(y,u(x))$ is
approximately differentiable at $x$ and
\begin{equation*}\label{eq:appdif}
\nabla v(x)=(\nabla_xg)(x,u)+(\nabla_rg)(x, u)\cdot\nabla u(x) \qquad
\text{$\L^{d}$-a.e. in $\R^d$\,;}
\end{equation*}

\item[(ii)] (jump part)
$J_v\subset J_u$
and it holds
\begin{equation*}\label{jump}
D^j v=\big(\widetilde g(x,u^+)-\widetilde g(x,u^-\big))\otimes\nu_{
u}\H^{d-1}\res J_u
\end{equation*}
in the sense of measures, where
$u^\pm(x)$ are the upper and lower approximate limits of
$u$ at $x$, {and $\widetilde g(x,r)$ denotes the precise representative of
$g(\cdot,r)$ on $\R^d\setminus \mathcal N$.}
\end{itemize}

Moreover
\begin{equation}\label{eq:appdif22}
\div\,v(x)\!=\!\left[(\div_xg)(x,u)+{\tr}\big((\nabla_rg)(x, u)\nabla u\big)\right]\!\L^{d}
+\langle\widetilde g(x,u^+)-\widetilde g(x,u^-),\nu_{
u}\rangle\H^{d-1}\res J_u
\end{equation}
in the sense of measures.
\end{thrm}


In order to obtain the analog in $GSBD$ of Theorem~\ref{chain rule}, we need a further assumption on $g$: 
\leavevmode
\begin{enumerate}[font={\normalfont},label={(G6)}]
\item for all
$x\in\R^d\setminus N$ the vector field $r\mapsto g(x,r)$ is \emph{conservative}; i.e., there exists a potential $G(x,\cdot)\in C^1(\R^d)$ such that $\nabla_r G(x,r)=g(x,r)$ for every $r\in\R^d$. \label{ass-f}
\end{enumerate}

\begin{thrm}\label{chain rule1}
Let $g$ be satisfying {\rm\ref{ass-a}}-{\rm\ref{ass-f}} above. Let $\Omega\subset\R^d$ be bounded. 
Then there exists a set $\mathcal N\subset\Omega$ with $\H^{d-1}(\mathcal N)=0$\,, such that,
 for every $r\in \R^d$ the function $g(\cdot,r)$ is approximately continuous in $\Omega\setminus \mathcal N$ and $\widetilde g(x,r)$ denotes the precise representative of
$g(\cdot,r)$ on $\Omega\setminus \mathcal N$, and
for any function $u\in GSBD^p(\Omega;\R^d)$,
the function $v(x):=g(x,u(x))$ is a {vector field whose distributional divergence is a Radon measure,} 
and the following nonautonomous chain rule formula for its measure divergence holds
%
%
%
%
%
%

\begin{equation}\label{eq:appdif22gsbd}
\div\,v(x)\!=\!\left[(\div_xg)(x,u)+(\nabla_rg)(x, u):e(u)\right]\!\L^{d}
+\langle\widetilde g(x,u^+)-\widetilde g(x,u^-),\nu_{
u}\rangle\H^{d-1}\res J_u
\end{equation}
in the sense of measures, i.e. for all open set $A\subseteq \Omega$  and for every $\varphi\in C^1_c(A)$
\begin{equation}
\begin{split}
-\int_A\langle g(x,u), \nabla \varphi\rangle\,\d x\! & =\!\int_A \varphi(\div_xg)(x,u)\,\d x+\int_A \varphi\big((\nabla_rg)(x, u):e(u)\big)\,\d x \\
& +\int_{A\cap J_u} \varphi\langle\widetilde g(x,u^+)-\widetilde g(x,u^-),\nu_{u}\rangle\, \d\H^{d-1}.
\end{split}
\label{eq:chainruleGSBD1}
\end{equation}
\end{thrm}
\begin{proof}
We may adapt the argument of \cite[Lemma~5.3]{FPS} to the nonautonomous setting as follows: we combine the nonautonomous formula proven in Theorem \ref{chain rule} for functions in $SBV$ with the approximation result in $GSBD^p$ recalled with Theorem~\ref{thm:cc}. 

Let $u\in GSBD^p(\Omega;\R^d)$ and $(u_k)\subset SBV^p(\Omega;\R^d)\cap L^\infty(\Omega;\R^d)$ be the sequence provided by Theorem~\ref{thm:cc}. Then, with \eqref{eq:appdif22}, 
for every $k\in \N$, for all open sets $A\subseteq \Omega$  and for every $\varphi\in C^1_c(A)$ we get
\begin{equation}
\begin{split}
-\int_A\langle g(x,u_k), \nabla \varphi\rangle\,\d x\!=&\!\int_A\varphi(\div_xg)(x,u_k)\,\d x+\int_A\varphi{\tr}\big((\nabla_rg)(x, u_k)\nabla u_k\big)\,\d x \\
& +\int_{A\cap J_{u_k}}\varphi\langle\widetilde g(x,u_k^+)-\widetilde g(x,u_k^-),\nu_{
u_k}\rangle\, \d\H^{d-1}.
\end{split}
\label{eq:chain1}
\end{equation}
Since $g(x,\cdot)$ is conservative, the matrix $\nabla_rg(x,\cdot)$ is symmetric. Then
$$
{\tr}\big((\nabla_rg)(x, u_k)\nabla u_k\big)= (\nabla_rg)(x, u_k):e(u_k) \,,
$$
and \eqref{eq:chain1} can be rewritten as
\begin{equation}
\begin{split}
-\int_A\langle g(x,u_k), \nabla \varphi\rangle\,\d x\!= & \!\int_A\varphi(\div_xg)(x,u_k)\,\d x+\int_A\varphi\big((\nabla_rg)(x, u_k):e(u_k)\big)\,\d x \\
&+\int_{A\cap J_{u_k}}\varphi\langle\widetilde g(x,u_k^+)-\widetilde g(x,u_k^-),\nu_{
u_k}\rangle\, \d\H^{d-1}\,.
\end{split}
\label{eq:chain2}
\end{equation}
Now, we aim to pass to the limit as $k\to+\infty$ in each of the four terms separately.

{As for the first term, by assumption {\rm\ref{ass-d}} we get
\begin{equation*}
|g(x, u_k)-g(x,u)|\leq M|u_k-u|\,,
\end{equation*}
whence $g(x, u_k)\to g(x,u)$ $\mathcal{L}^d$-a.e. in $A$. Then, with {\rm\ref{ass-a}} 
and the boundedness of $\|\nabla\varphi\|_\infty$ we infer
\begin{equation}
\lim_{k\to+\infty} \int_A\langle g(x,u_k), \nabla \varphi\rangle\,\d x = \int_A\langle g(x,u), \nabla \varphi\rangle\,\d x
\label{eq:5.7}
\end{equation}
by the dominated convergence theorem.} To prove that
\begin{equation}
\lim_{k\to+\infty} \int_A\varphi(\div_xg)(x,u_k)\,\d x = \int_A\varphi(\div_xg)(x,u)\,\d x
\label{eq:5.7bis}
\end{equation}
we notice that by assumption {\rm\ref{ass-b}}, $(\div_xg)(x,u_k)\to(\div_xg)(x,u)$ $\mathcal{L}^d$-a.e. and the sequence $((\div_xg)(x,u_k)-(\div_xg)(x,u))_k$ is dominated by an $L^1$-function. Then \eqref{eq:5.7bis} follows again from the dominated convergence theorem since $\varphi$ is bounded.

By assumptions {\rm\ref{ass-d}} and {\rm\ref{ass-e}}, Theorem~\ref{thm:cc}(i) and the dominated convergence theorem we infer that $(\nabla_rg)(\cdot, u_k)\to(\nabla_rg)(\cdot, u)$ in $L^q$ for any $q\in[1,+\infty)$. Then, with Theorem~\ref{thm:cc}(ii), $\varphi\in C^1_c(A)$ and H\"older's inequality we finally get
\begin{equation}
\lim_{k\to+\infty} \int_A\varphi\big((\nabla_rg)(x, u_k):e(u_k)\big)\,\d x = \int_A\varphi\big((\nabla_rg)(x, u):e(u)\big)\,\d x\,. 
\label{eq:5.8}
\end{equation}
Finally, since $r\to \widetilde g(x,r)$ is Lipschitz continuous (see \cite[Proposition~3.2(i)]{ACDD}), the proof of
\begin{equation}
\lim_{k\to+\infty} \int_{A\cap J_{u_k}}\varphi\langle\widetilde g(x,u_k^+)-\widetilde g(x,u_k^-),\nu_{
u_k}\rangle\, \d\H^{d-1} = \int_{A\cap J_{u}}\varphi\langle\widetilde g(x,u^+)-\widetilde g(x,u^-),\nu_{
u}\rangle\, \d\H^{d-1}
\label{eq:5.9}
\end{equation}
is exactly as in \cite{FPS}, by exploiting Theorem~\ref{thm:cc} $(iv)$, so we omit the details.
\end{proof}

\section{Nonautonomous symmetric jointly convex functions}\label{sec: joint}

We give a definition of nonautonomous (NA) symmetric jointly convex function with $W^{1,1}$ dependence of the approximating vector fields with respect to the spatial variable $x$\,. This can be considered as an extension to the nonautonomous setting of the definition of {\it symmetric jointly convex function} (see Definition~\ref{defn:sjconv} below), recently introduced in \cite{FPS}.

{
\begin{dfntn}[Symmetric joint convexity]\label{defn:sjconv}
Let  $\phi: \R^d\times \R^d\times \R^{d}\to  [0,+\infty)$. We say that
$\phi$ is {\it symmetric jointly convex} if there exists a sequence of uniformly continuous, bounded, conservative vector fields $g_j:\R^d\to \R^d$ such that 
\begin{equation*}
\phi(r,t,\xi)=\sup_{j\in\N} \langle g_j(r)-g_j(t),\xi\rangle\quad\quad
{\rm for\ all\ } (r,t,\xi)\in \R^d\times \R^d\times  \R^{d}\,.
\end{equation*}
\end{dfntn}
The definition of nonautonomous (NA) symmetric jointly convex function is the following.}

\begin{dfntn}[NA symmetric joint convexity]\label{defn:nasjconv}
Let  $\phi: \Omega\times \R^d\times \R^d\times \R^{d}\to  [0,+\infty)$. We say that
$\phi$ is {\it NA symmetric jointly convex} if there exists a sequence of bounded functions $g_j:\Omega\times  \R^d\to \R^d$ such that 
\begin{equation}
\phi(x,r,t,\xi)=\sup_{j\in\N} \langle g_j(x,r)-g_j(x,t),\xi\rangle\quad\quad
{\rm for\ all\ } (x,r,t,\xi)\in \Omega\times \R^d\times \R^d\times  \R^{d}
\label{eq:reprjointlyconv}
\end{equation}
and for every $j\in\N$  the function $g_j$ satisfies conditions {\rm\ref{ass-a}}, {\rm\ref{ass-b}}, {\rm\ref{ass-f}} and the following condition
\leavevmode

\begin{enumerate}[font={\normalfont},label={(G3$^\prime$)}]
\item 
there exists a Lipschitz constant $L_j$ independent of $x$ such that 
\[
 |g_j(x,r)-g_j(x,t)|\leq L_j|r-t|
 \]
for all $r,\,t\in \R^d$ and {$\mathcal{H}^{d-1}$-a.e. $x\in\R^d$}\,. \label{ass-c'}
\end{enumerate}
\end{dfntn}

\begin{rmrk}\label{precis}
An inspection {of} the proof of \cite[Lemma~2.4]{DCFV1} shows that, under assumptions {\rm\ref{ass-a}}, {\rm\ref{ass-b}} and {\rm\ref{ass-c'}}, for every $j\in \N$ there exists a subset $N_j\subset\Omega$ (indipendent of $r$) with $\H^{d-1}(N_j)=0$, such that for any $x\in\Omega\backslash N_j$ and any $r\in\R^d$, $g_j(\cdot,r)$ is approximately continuous at $x$. This implies, in particular, that $\widetilde g_j(x,r)=g_j(x,r)$ for every $x\in\Omega\setminus N_j$ and $r\in \R^d$\,. Let $N=\cup_{j\in\N} N_j$, then $\H^{d-1}(N)=0$. In the following we will tacitly assume that the integrand $\phi(\cdot,r,t,\xi)$ coincides in $\Omega\setminus N$ (and so $\H^{d-1}$-a.e. in $\Omega$) with the following representative 
$$
\overline\phi(\cdot,r,t,\xi):=\sup_{j\in\N} \langle \widetilde g_j(\cdot,r)-\widetilde g_j(\cdot,t),\xi\rangle\quad\quad
{\rm for\ all\ } (r,t,\xi)\in \R^d\times \R^d\times  \R^{d}.
$$
\end{rmrk}

\subsection{Some examples}
\label{r:r14}

We give some example of NA symmetric jointly convex functions (see also \cite[Remark~3.2]{DC}).
A first example is the model case, i.e.,
\begin{equation}\label{A0}
 \phi(x,r,t,\xi):=\langle g(x,r)-g(x,t),\xi\rangle^+\,, \tag{$\mathcal{A}0$}
\end{equation}
where $s^+:=\max\{s,0\}$ for $s\in\R$, and $g$ satisfies conditions {\rm\ref{ass-a}}, {\rm\ref{ass-b}}, {\rm\ref{ass-c'}}  and {\rm\ref{ass-f}}.
A further example is
\begin{equation}\label{A1}
\phi(x,r,t,\xi) :=a(x)\kappa(r,t,\xi)\,, \tag{$\mathcal{A}1$}
\end{equation}
where $a$ is a nonnegative bounded $W^{1,1}$ function, $$
\kappa(r,t,\xi)=\sup_{j\in\N}\langle h_j(r)-h_j(t),\xi\rangle^+
$$and $h_j$ is a sequence of continuous functions and for every $j\in\N$ the function $h_j$ is a conservative vector field according to {\rm\ref{ass-f}}.

With the examples of (autonomous) symmetric jointly convex integrands $\phi(r,t,\xi)$ present in literature at hand, we can immediately construct examples of nonautonomous symmetric jointly convex integrands:
\begin{itemize}
\item[$(i)$] $\phi(x,r,t,\xi) :=a(x)\kappa(|\langle r-t, \xi \rangle|)$, where $\kappa$ is a non-negative convex, subadditive and increasing function (see \cite{GZ1} for the case $a(x)\equiv 1$);
\item[$(ii)$] $\displaystyle\phi(x,r,t,\xi) :=a(x)\sup_{({\zeta}_1,\dots,{\zeta}_k)}\left(\sum_{k=1}\theta_k(\langle r-t, {\zeta}_k\rangle)^2|\langle \xi, {\zeta}_k\rangle|^2\right)^{1/2}$, where $\theta_k$, $k=1,\dots,d$, are even, continuous, subadditive functions such that $\theta_k(0)=0$, and the supremum is taken over all orthonormal bases $({\zeta}_k)_{k=1}^d$ of $\R^d$. The symmetric joint convexity of the corresponding $\kappa(r,t,\xi)$ has been proved in \cite[Proposition~4.5]{FPS} (see also \cite{DOT}); 
\item[$(iii)$] $\phi(x,r,t,\xi) :=a(x)\gamma(|r-t|)|\xi|$, which can be seen of type \ref{A1} with $\kappa(r,t,\xi):=\gamma(|r-t|)|\xi|$ (see \cite[Theorem~4.1]{FPS}); 
\item[$(iv)$] $\phi(x,r,t,\xi):= a(x)\kappa(\xi)$, where $\kappa$ is even, 1-homogeneous and convex (see \cite[Proposition~4.13]{FPS}). 
\end{itemize}

The following result deals with more general nonautonomous integrands independent of the traces $r,t$; i.e.,
\begin{equation}
\phi(x,r,t,\xi) = \kappa(x,\xi) \qquad \mbox{ for all } (x,r,t,\xi)\in\Omega\times\R^d\times\R^d\times\R^d\,,\,\, \mbox{ with $r\neq t$.}
\label{eq:jointconv4.13}
\end{equation}
It can be seen as a further generalization of \cite[Proposition~4.13]{FPS}.

\begin{prpstn}\label{thm:thm4.13fps}
Let $\phi:\Omega\times\R^d\times\R^d\times \R^d\to[c,+\infty)$, $c>0$, be of type \eqref{eq:jointconv4.13}, where $\kappa:\Omega\times\R^d\to[0,+\infty)$ 
complies with 
\begin{enumerate}[font={\normalfont},label={(K\arabic*)}]
\item $x\mapsto \kappa(x,\xi)$ belongs to $W^{1,1}(\Omega;\R^d)$ for all $\xi\in\R^d$\,; \label{ass-k1}
\item there exist a positive function $h\in L^1_{\rm loc}(\R^d)$ and a modulus of continuity $\omega$ such that
\[
|\nabla_x \kappa(x,\xi)-\nabla_x \kappa(x,\xi')|\le
\omega(|\xi-\xi'|)h(x)
\]
for all $\xi,\,\xi'\in \R^d$ and  for $\L^{d}$-a.e. $x\in\Omega$\,; \label{ass-k2}
\item $\kappa(x,\cdot)$ is even, positively 1-homogeneous and convex for every $x\in\Omega$. \label{ass-k3}
\end{enumerate}
Then 
$\phi$ is {NA} symmetric jointly convex.
\end{prpstn}
\begin{proof}

First of all, we observe that since $\kappa$ is locally bounded in $\Om\times\R^d$
and positively $1$-homogeneous with respect to $\xi$, for any open set
$U\subset\!\subset\Om$, there exists a constant $\Lambda_U$
such that
\begin{equation}\label{e:eqnumero1}
\qquad\qquad
0\leq \kappa(x,\xi)\leq\Lambda_U|\xi| \qquad \hbox{for all $(x,\xi)\in U\times\R^d\,.$}
\end{equation}
Then, the convexity of $\kappa$ with respect to $\xi$ combined with \eqref{e:eqnumero1} yields  that
\begin{equation}\label{semi22bis}
\qquad\qquad
|\kappa(x,\xi_1)-\kappa(x,\xi_2)|\leq \tilde{c}\Lambda_U|\xi_1-\xi_2| \qquad
\hbox{for all $(x,\xi_1),(x,\xi_2)\in U\times\R^d\,,$}
\end{equation}
for some dimensional constant $\tilde{c}>0$. Now, let us fix a dense sequence $\{\xi_n\}\subset\R^d$.
Thanks to \ref{ass-k1} for all $n\in \N$ there exists a {Borel} set $V_{n}\subset U$, with $\H^{d-1}(V_{n})=0$,
such that $\kappa(\cdot,\xi_n)$ is approximately continuous in $\Om\setminus V_{n}$.
Setting $V=\cup_{n}V_{n}$ we obtain that $\H^{d-1}(V)=0$.
Making use of \eqref{semi22bis},
one easily gets that $\kappa(\cdot,\xi)$ is approximately  continuous  in
$U\setminus V$ for all $\xi\in\R^d$. Then $\kappa(\cdot,\xi)=\widetilde \kappa(\cdot,\xi)$ for $\H^{d-1}$-a.e. $x\in\Om$ and for all $\xi\in\R^d$.
Similarly, by using \ref{ass-k2} there exists
a {Borel} set $M\subset U$, with $\L^{d}(M)=0$,
such that  $\kappa(\cdot,\xi)$ is approximately  differentiable  in
$U\setminus M$ for all $\xi\in\R^d$.


Now we may adapt the argument of \cite[Proposition~4.13]{FPS} to the nonautonomous setting by using the De Giorgi's approximation result, Lemma~\ref{degiorgi}. Indeed, by virtue of this result, 
for $\H^{d-1}$-a.e. $x\in\Omega$ and for every $\xi\in\R^d$ we can write
\begin{equation}
\kappa(x,\xi)=\sup_{j\in\N}\langle a_{j}(x),\xi\rangle\,,
\label{eq:degiorgirepr}
\end{equation}
where the functions $a_{j}$ are defined as in \eqref{coeff1} with $f=\kappa$. 
Notice that the functions $a_{j}$ are bounded and $a_{j}\in W^{1,1}(\Omega;\R^d)$. Moreover, as proven in \cite[Lemma~2.4]{DCFV1}, we have
$$
a_{j}(x)=\widetilde{a}_{j}(x)=-\int_{\R^{d}}\widetilde \kappa(x,z)\nabla\alpha_{j}(z)\,\mathrm{d}z
$$
and
\begin{equation}
\kappa(x,\xi)=\sup_{j\in\N}\langle {\widetilde a}_{j}(x),\xi\rangle\,,
\label{eq:degiorgireprtilde}
\end{equation}
for $\H^{d-1}$-a.e. $x\in\Om$ and for all $\xi\in\R^d$.
Moreover, for every $i=1,\dots,d$
$$
\partial_{x_i} a_{j} (x)=-\int_{\R^{d}}\partial_{x_i} \kappa(x,z)\nabla\alpha_{j}(z)\,\mathrm{d}z
$$
for $\L^{d}$-a.e. $x\in\Om$. 

Let $N\subset\Omega$ be the set such that $\mathcal{H}^{d-1}(N)=0$ and \eqref{eq:degiorgireprtilde} holds for $x\in\Omega\backslash N$. Correspondingly, for every such $x$ we define the countable set $\mathcal{A}(x)$ as in \eqref{eq:classdomain}.

For every $h, l\in \N$, $b_l:{\Omega\backslash N}\to\R^d$ any selection of $\mathcal{A}$ and $p\in\R^d$, we set
\begin{equation}
g_{h,l,p}(x,w):=\theta_h(\langle w-p, b_{l}(x)\rangle) b_{l}(x)\,,\quad w\in\R^d\,,\, x\in\Omega\backslash N\,,
\label{eq:w11vectorfiels}
\end{equation}
where $\theta_h:\R\to[0,+\infty)$ is {defined as $\theta_h(y):=\frac{2}{\pi}\arctan (h |y|)$. Note that $\theta_h$ is a $C^1(\R\backslash\{0\})$, subadditive function such that
$\theta_h(0)=0$, $\theta_h\leq1$ and $\theta_h(y)\to 1$ as $h\to+\infty$ for each fixed $y\neq0$.} 
Clearly, each $g_{h,l,p}$ is bounded and Lipschitz continuous in $w$ uniformly with respect to $x$, 
so that, taking into account also the remarks about $a_{j}$ above, conditions {\rm\ref{ass-a}}, {\rm\ref{ass-b}} and {\rm\ref{ass-c'}} are satisfied. Moreover, $g_{h,l,p}$ is conservative with potential
\begin{equation*}
G_{h,l,p}(x,w):=\Theta_h\big( \langle w-p,b_{l}(x)\rangle \big)\,,\quad w\in\R^d\,,\,x\in\Omega\backslash N\,,
\end{equation*}
where $\Theta_h$ is a primitive of $\theta_h$. Thus, also {\rm\ref{ass-f}} is satisfied. 

We claim that for every $x\in\Om\backslash N$ and for all $(r,t,\xi) \in \R^d\times \R^d\times \R^d$ with $r \neq t$ we have
\begin{equation}
\phi(x,r,t,\xi)=\kappa(x,\xi) =  \sup_{h,l,p} \langle g_{h,l,p}(x,r)-g_{h,l,p}(x,t) , \xi \rangle\,.
\label{eq:representationformul}
\end{equation}
We first prove the inequality ``$\ge$'' in \eqref{eq:representationformul}. Let $x\in\Omega\backslash N$ be fixed. For each   $(r,t,\xi) \in \R^d\times \R^d\times \R^d$ with $r \neq t$, we get from the definition of $b_l$, the subadditivity of $\theta_h$ and $\theta_h\leq1$ that
\begin{align*}
\left\langle g_{h,l,p}(x,r)- g_{h,l,p}(x,t),\xi  \right\rangle \leq \langle b_{l}(x),\xi\rangle\leq \kappa(x,\xi) = \phi(x,r,t,\xi)\,,
\end{align*}
whence the desired inequality follows by passing to the supremum on the left hand side.

Now, we turn to the proof of the reverse inequality in \eqref{eq:representationformul}. 
Let $(p_m)_{m\in\N}\subset \mathbb{Q}^d$ be a sequence such that $p_m\to t$ as $m\to+\infty$. Then, with fixed $h,l$, with the continuity of $\theta_h$ and recalling that $\theta_h(0)=0$, we have
\begin{equation}
\begin{split}
\lim_{m\to+\infty} \left\langle g_{h,l,p_m}(x,r)- g_{h,l,p_m}(x,t),\xi  \right\rangle & = \left\langle \left(\theta_h(\langle r-t, b_{l}(x)\rangle)-\theta_h(0)\right)b_{l}(x), \xi \right\rangle \\
 & = \theta_h(\langle r-t, b_{l}(x)\rangle)\langle b_{l}(x),\xi\rangle\,. 
\end{split}
\label{eq:stimacontinuap}
\end{equation}
Let $\eps>0$ be fixed. {We claim that there exists $l_0\in\N$ such that 
\begin{equation}
\langle b_{l_0}(x),\xi\rangle > \kappa(x,\xi)-\eps \quad \mbox{ and } \quad |\langle b_{l_0}(x), r-t \rangle|\neq0\,.
\label{eq:stimacontinua2bis}
\end{equation}
 First, by \eqref{eq:degiorgirepr} there exists $j_0$ such that $\langle a_{j_0}(x),\xi\rangle > \kappa(x,\xi)-\eps$. We may assume also that $|\langle a_{j_0}(x), r-t \rangle|\neq0$ and set $b_{l_0}(x):=a_{j_0}(x)$. If not, we may replace $a_{j_0}(x)$ by $b_{l_0}(x)\in \mathcal{A}(x)$ defined as
\begin{equation*}
b_{l_0}(x):= a_{j_0}(x) + \sigma_{j_0}{\rm sign}(\langle v_{j_0}, \xi\rangle) v_{j_0}
\end{equation*}
for some $v_{j_0}\in\mathbb{Q}^d\cap\mathbb{S}^{d-1}$ such that $\langle v_{j_0},r-t\rangle\neq0$, $\langle v_{j_0}, \xi\rangle\neq0$ and some $\sigma_{j_0}\in\mathbb{Q}^+$. We then have
\begin{equation*}
\langle b_{l_0}(x), r-t\rangle = \langle a_{j_0}(x), r-t \rangle + \sigma_{j_0}{\rm sign}(\langle v_{j_0}, \xi\rangle) \langle v_{j_0}, r-t\rangle = \pm \sigma_{j_0} \langle v_{j_0}, r-t\rangle \neq 0\,,
\end{equation*}
and
\begin{equation*}
\langle b_{l_0}(x), \xi \rangle = \langle a_{j_0}(x) , \xi \rangle + \sigma_{j_0}{\rm sign}(\langle v_{j_0}, \xi\rangle) \langle v_{j_0}, \xi\rangle > \langle a_{j_0}(x),\xi\rangle > \kappa(x,\xi)-\eps\,.
\end{equation*}
This proves \eqref{eq:stimacontinua2bis}. Now, since $\theta_h(y)\to1$ as $h\to+\infty$ for each fixed $y\neq0$,
\begin{equation}
\lim_{h\to+\infty} \theta_h(\langle r-t, b_{l_0}(x)\rangle)\langle b_{l_0}(x),\xi\rangle = \langle b_{l_0}(x),\xi\rangle\,.
\label{eq:stimacontinua2}
\end{equation}
Combining the previous estimates \eqref{eq:stimacontinuap}, \eqref{eq:stimacontinua2bis} and \eqref{eq:stimacontinua2}, we then get
\begin{equation*}
\begin{split}
\sup_{h,l,p} \langle g_{h,l,p}(x,r)-g_{h,l,p}(x,t) , \xi \rangle & \geq \lim_{h\to+\infty}  \theta_h(\langle r-t, b_{l_0}(x)\rangle)\langle b_{l_0}(x),\xi\rangle \\
& >\kappa(x,\xi)-\eps\,.
\end{split}
\end{equation*}}
The arbitrariness of $\eps$ yields ``$\le$'' in \eqref{eq:representationformul}, and this concludes the proof. 
\end{proof}

\subsection{Lower semicontinuity result}

We are in position to prove the first main result of our paper, Theorem~\ref{thm:gsbdsemicontinuity1} below: the $NA$ symmetric joint convexity of the integrand, introduced with Definition~\ref{defn:nasjconv}, is a sufficient condition for the lower semicontinuity in $GSBD^p$ of the corresponding functional. 
We will first show the assertion for $NA$ symmetric jointly convex $\phi$ admitting a representation \eqref{eq:reprjointlyconv} with vector fields $g_j$ satisfying the assumptions of the nonautonomous chain rule for its divergence, Theorem~\ref{chain rule1}. This will allow us to follow some arguments contained in the proofs of \cite[Theorem~3.4]{DC} and \cite[Theorem~5.1]{FPS}.  

The general case of $g_j$ complying with assumptions {\rm\ref{ass-a}}, {\rm\ref{ass-b}}, {\rm\ref{ass-c'}} and {\rm\ref{ass-f}}, which are the only needed in Definition~\ref{defn:nasjconv}, then will follow through a standard mollification technique (see, e.g., \cite[Theorem~5.22]{AFP}).
\begin{thrm}\label{thm:gsbdsemicontinuity1}
Let  $\phi: \Omega\times \R^d\times \R^d\times \R^{d}\to  [0,+\infty)$ be a 
 NA symmetric jointly convex function.
Then, for every $(u_n)_n\subset GSBD^p(\Om;\R^d)$, $p>1$, converging in measure to $u\in GSBD^p(\Om;\R^d)$ 
and such that
\begin{equation}
C:=\sup_{n\in\N}\left[\int_\Om |e(u_n)|^p\,\mathrm{d}x + \H^{d-1}(J_{u_n})\right] <+\infty\,,
\label{eq:equibounded}
\end{equation}
we have
\begin{equation}\label{eq:e17bbis}
\int_{J_u\cap\Om} \phi(x,u^+,u^-,\nu_{u})\,\mathrm{d}\H^{d-1} \leq
\liminf_{n\to +\infty}\int_{J_{u_n}\cap\Om} \phi(x,u_n^+,u_n^-,\nu_{u_n})\,\mathrm{d}\H^{d-1}\,.
\end{equation}
\end{thrm}

\begin{proof}
Assume that \eqref{eq:reprjointlyconv} holds for some sequence of vector fields $(g_j)_{j\in\N}$. We subdivide the proof into two steps.

\emph{Step~1: each $g_j$ satisfies conditions {\rm\ref{ass-a}} - {\rm\ref{ass-f}}.} 
Let $(u_n)_n$ and $u$ be as above. Since $\phi$ is non-negative, by virtue of Lemma~\ref{l:passup} it will suffice to show that
\begin{equation}
\liminf_{n\to+\infty}F_{g}(u_n)\geq
F_{g}(u)\,,
\label{eq:5.11}
\end{equation}
where
$$
F_{g}(w):=\int_{J_w\cap A} \langle g(x,w^+)-g(x,w^-),\nu_w\rangle^+\,\mathrm{d}\mathcal{H}^{d-1}\,,
$$
for every vector field $g$ complying with {\rm\ref{ass-a}}-{\rm\ref{ass-f}} and every open set $A\subset\Omega$. 
{ 
Since 
$$
F_{g}(u)=\sup\left\{ \int_{J_u\cap A} \langle g(x,u^+)-g(x,u^-),\nu_u\rangle\psi\, \mathrm{d}\H^{d-1}: \,\, 
\psi\in C^1_c(A),\ 0\leq\psi\leq 1\right\},
$$
again by virtue of Lemma~\ref{l:passup} the lower semicontinuity of the functional $F_{g}$ will follow if we prove the continuity of
\begin{equation}\label{p:p113}
F_g^\psi(u):= \int_{J_u} \langle g(x,u^+)-g(x,u^-),\nu_u\rangle\psi\, \mathrm{d}\H^{d-1}
\end{equation}
for some fixed $\psi\in C^1_0(A)$, $0\leq\psi\leq1$.
Indeed, by using the chain rule formula \eqref{eq:chainruleGSBD1} 
we have
\begin{equation*}
\begin{split}
\int_{A\cap J_u}\, \langle g(x,u^+)-g(x,u^-),\nu_u\rangle\,\psi\, \mathrm{d}\H^{d-1}= & - \int_A\psi(\div_xg)(x,u)\,\d x - \int_A\psi\big((\nabla_rg)(x,u):e(u)\big)\,\d x \\
& - \int_A\langle g(x,u), \nabla \psi\rangle\,\d x\,.
\end{split}
\end{equation*}
We claim that
\begin{alignat}2
\label{eq:e30}
& \int_A\langle g(x,u), \nabla \psi\rangle\,\d x
= \lim_{n\to +\infty} \int_A\langle g(x,u_n), \nabla \psi\rangle\,\d x\,,
\\
\label{eq:e31}
& \int_A\psi(\div_xg)(x,u)\,\d x
= \lim_{n\to +\infty} \int_A\psi(\div_xg)(x,u_n)\,\d x\,,
\\
\label{eq:e32}
& \int_A\psi\big((\nabla_rg)(x,u):e(u)\big)\,\d x
= \lim_{n\to +\infty} \int_A\psi\big((\nabla_rg)(x,u_n):e(u_n)\big)\,\d x\,.
\end{alignat}
As for \eqref{eq:e30}, from assumptions {\rm\ref{ass-a}} and {\rm\ref{ass-d}} we get $g(x, u_n)\to g(x,u)$ for $\mathcal{L}^d$-a.e. $x$ in $A$, and
\begin{equation*}
|\nabla\psi||g(x,u_n)|\leq \|\nabla\psi\|_\infty h_1(x) \quad \mbox{ for $\mathcal{L}^d$-a.e. $x$ in $A$, for every $n\in\N$,}
\end{equation*}
whence the assertion follows from the dominated convergence theorem. 
The proof of \eqref{eq:e31} is similar to that of \eqref{eq:5.7bis}. In order to prove \eqref{eq:e32}, we first observe that under assumption \eqref{eq:equibounded}, we can find a subsequence (not relabeled) $(u_n)$ such that $e(u_n)\rightharpoonup e(u)$ weakly in $L^1(A;\mathbb M_{\mathrm{sym}}^{d\times d})$. Then, by using {\rm\ref{ass-d}}-{\rm\ref{ass-e}} and dominated convergence, we have $\nabla_{r} g(x,u_n)\to \nabla_{r} g(x,u)$ strongly in $L^{1}(A;\mathbb M^{d\times d})$. Now, writing
\begin{equation*}
(\nabla_rg)(x,u_n) : e(u_n) = [(\nabla_rg)(x,u_n)-(\nabla_rg)(x,u)]:e(u_n) + (\nabla_rg)(x,u) : e(u_n)\,,
\end{equation*}
the first term in the right hand side tends to 0 in $L^{1}(A;\mathbb M^{d\times d})$ by virtue of dominated convergence theorem and {\rm\ref{ass-d}}. 
This implies that
\begin{equation*}
(\nabla_rg)(x,u_n):e(u_n) \rightharpoonup (\nabla_rg)(x,u):e(u) \quad \mbox{ in $L^{1}(A;\mathbb M^{d\times d})$, }
\end{equation*}
and since $\psi\in L^\infty(A)$, \eqref{eq:e32} follows.
Collecting \eqref{eq:e30}, \eqref{eq:e31} and \eqref{eq:e32}, we conclude that the functionals $F^\psi_g$, $\psi\in C^1_c(A)$ are continuous along the sequence $(u_n)$, and so the lower semicontinuity of $F_g$ follows\,.
}


\emph{Step~2: each $g_j$ satisfies {\rm\ref{ass-a}}, {\rm\ref{ass-b}}, {\rm\ref{ass-c'}}, and {\rm\ref{ass-f}}.} 
By Lemma~\ref{l:passup} again, it is sufficient to prove \eqref{eq:5.11}
for every vector field $g$ complying with {\rm\ref{ass-a}}, {\rm\ref{ass-b}}, {\rm\ref{ass-c'}} and {\rm\ref{ass-f}}, and every open set $A\subset\Omega$. 

For this, given a non-negative, even function $\rho\in C^\infty_c(\R^d)$ such that ${\rm supp}\,\rho\subset B_1$ and $\int_{\R^d}\rho\,\mathrm{d}t=1$, we define a sequence of mollifiers $(\rho_\eps)_\eps$ by setting $\rho_\eps(t):=\frac{1}{\eps^d}\rho(\frac{t}{\eps})$. Correspondingly, we consider the mollified functions
$$
g_\varepsilon(x,r):= (g*\rho_\eps)(x,r)=\int_{\R^d}g(x,r-t)\rho_\varepsilon(t)\mathrm{d}t\,.
$$
Then each $g_\varepsilon$ satisfies {\rm\ref{ass-a}}--{\rm\ref{ass-f}}. Indeed, 
for all
$r\in\R^d$ the function $x\mapsto g_\varepsilon(x,r)$ belongs to $W^{1,1}_{\rm loc}(\R^d;\R^d)$, we have $\nabla_xg_\varepsilon(x,r) = (\nabla_xg*\rho_\eps)(x,r)$  for all $r\in \R^d$ and $\mathcal{L}^d$-a.e. $x\in\R^d$,
and
$$
|\nabla_x g_\varepsilon(x,r)-\nabla_x g_\varepsilon(x,s)|\le
\int_{\R^d}|\nabla_x g(x,r-t)-\nabla_x g(x,s-t)|\rho_\varepsilon(t)\,\mathrm{d}t\le
h(x)\omega(|r-s|)
$$
for all $r, s\in\R^d$ and $\mathcal{L}^d$-a.e. $x\in\R^d$. Moreover, the function $r\mapsto g_\varepsilon(x,r)$ is continuously differentiable in $
\R^d$ for $\mathcal{H}^{d-1}$-a.e. $x\in\R^d$,
 $$|\nabla_r g_\varepsilon(x,r)|=|(g*\nabla\rho_\eps)(x,r)|\leq M\,,$$ 
and by {\rm\ref{ass-c'}}
\[
 |\nabla_r g_\varepsilon(x,r)-\nabla_r g_\varepsilon(x,s)|\leq C\tilde \omega(|r-s|)
 \]
for all $r,\,s\in \R^d$ and $\mathcal{H}^{d-1}$-a.e. $x\in\R^d$\,.
In order to prove that $g_\varepsilon(x,r)$ is a conservative field it suffices to define its potential in the following way
$$
G_\varepsilon(x,r):=(G*\rho_\varepsilon)(x,r)=\int_{\R^d}G(x,t)\rho_\varepsilon(r-t)\,\mathrm{d}t
$$
where for $\mathcal{H}^{d-1}$-a.e. $x\in\R^d$, $G(x,\cdot)\in C^1(\R^d)$ is a potential of $g(x,r)$ for every $r\in\R^d$. 
Indeed,
$$
\nabla_r G_\varepsilon(x,r)=\int_{\R^d}\nabla_rG(x,r-t)\rho_\varepsilon(t)\mathrm{d}t=\int_{\R^d}g(x,r-t)\rho_\varepsilon(t)\mathrm{d}t=g_\varepsilon(x,r).
$$
Therefore, from Step~1 we can deduce the lower-semicontinuity result \eqref{eq:5.11} for each $g_\eps$, namely
\begin{equation}
\liminf_{n\to+\infty}F_{g_\varepsilon}(u_n)\geq F_{g_\varepsilon}(u)
\label{eq:5.11bis}
\end{equation}
for all $\varepsilon>0$. Moreover, by {\rm\ref{ass-c'}}, we have $|g(x,r)-g_\varepsilon(x,r)|\leq \hat \omega(\varepsilon)$ and so 
by \eqref{eq:5.11bis} and \eqref{eq:equibounded} we get
$$
\liminf_{n\to+\infty}F_{g}(u_n)\geq
\liminf_{n\to+\infty}[F_{g_\varepsilon}(u_n)-2\hat \omega(\varepsilon)\H^{d-1}(J_{u_n})]\geq
F_{g_\varepsilon}(u)-2\hat \omega(\varepsilon)C.
$$
The desired assertion then follows letting $\varepsilon\to 0$, and this concludes the proof.
\end{proof}

\section{Nonautonomous $ {{BV}}$ symmetric jointly convex functions}\label{sec: bvjoint}

In this section we give a new definition of nonautonomous symmetric joint convexity for functions with $BV$ dependence with respect to the spatial variable $x$\,. This definition is new 
even for functions with $W^{1,1}$ dependence. Indeed, as it will be apparent from Definition~\ref{defn:bvsjconv} below, the $W^{1,1}$ regularity in $x$ of the integrand $\phi$ is assumed, while in Section~\ref{sec: joint} 
the analogous assumption is required for the approximating vector fields.
The new definition
coincides with the $NA$ symmetric joint convexity for splitting-type integrands of the form \eqref{A0} or \eqref{A1}. Moreover, the lower semicontinuity theorem for $NA$ symmetric jointly convex functions (Theorem~\ref{thm:gsbdsemicontinuity1}) will be a key tool for the proof of the Proposition~\ref{corcor} which is a first step in order to obtain the analogous result for nonautonomous BV symmetric jointly convex functions, see Theorem~\ref{mmmdd}.

\begin{dfntn} \label{defn:bvsjconv}
A function $\phi: \Omega\times \R^d\times \R^d\times \R^{d}\to  [0,+\infty)$ is said to be
$BV$  {\it symmetric jointly convex} 
if the following conditions hold:
\leavevmode
\begin{enumerate}[font={\normalfont},label={(B{\arabic*})}]
\item for every $(r,t,\xi)\in \R^d\times \R^d\times \R^{d}$ the function $\phi(\cdot,r,t,\xi)$ belongs to $BV$ and there exists a Borel set $N\subset \Omega$ with $\H^{d-1}(N)=0$ such that $\phi(\cdot,r,t,\xi)$ coincides with its lower approximate limit $\phi^-(\cdot,r,t,\xi)$ in $\Omega\setminus N$ for all $(r,t,\xi)\in \R^d\times \R^d\times \R^{d}$; \label{ass-b1}
\item for every $x\in\Omega\setminus N$ the function $\phi(x,\cdot,\cdot,\cdot)$ is symmetric jointly convex\,; \label{ass-b2}
\item there exists $L>0$ such that
\[
 |\phi(x,r,t,\xi)-\phi(x,s,t,\xi)|\leq L|r-s|
 \] 
 for all $x\in\Omega\setminus N$, for all $r,s,t,\xi\in \R^d$. \label{ass-b3}
\end{enumerate}
\end{dfntn}

As we will prove with Theorem~\ref{mmmdd} at the end of this section, a lower semicontinuity result in $GSBD^p$ for ${{BV}}$ symmetric jointly convex integrands can be obtained by requiring the further condition that $\phi$ be strictly positive for $\H^{d-1}$-a.e. $x\in\Omega$.

\subsection{Some examples}
\label{r:r15}

The examples of Section~\ref{r:r14} can be easily adapted in order to construct $BV$ symmetric jointly convex functions.

A first example is the model case \eqref{A0}, 
with $g$ satisfying \ref{ass-c'} 
and \ref{ass-f} and the following condition:
\begin{enumerate}[font={\normalfont},label={(G1$^\prime$)}]
\item for every $r\in \R^d$ the function $g(\cdot,r)$ is a locally bounded $BV$ and there exists a Borel set $N\subset \Omega$ with $\H^{d-1}(N)=0$ such that $g(\cdot,r)=g^-(\cdot,r)$ in $\Omega\setminus N$ for all $r\in \R^d$\,. \label{ass-C1}
\end{enumerate}


Another example is \eqref{A1}
where $a$ is a nonnegative, bounded $BV$ function coinciding with its lower approximate limit $a^-$, and $\kappa$ is a symmetric jointly convex function {(see Definition~\ref{defn:sjconv})} for which there exists $L>0$ such that
\[
 |\kappa(r,t,\xi)-\kappa(s,t,\xi)|\leq L|r-s|
 \] 
for all $r,s,t\in \R^d$ and $\xi\in\R^d$.

The following is the $BV$ counterpart of Proposition~\ref{thm:thm4.13fps}.

\begin{prpstn}\label{thm:thm4.13fpsbv}
Let $\phi:\Omega\times\R^d\times\R^d\times\R^d\to[c,+\infty)$, $c>0$, be of type \eqref{eq:jointconv4.13}, where $\kappa:\Omega\times\R^d\to[0,+\infty)$ 
complies with {\rm\ref{ass-k3}} and
\begin{enumerate}[font={\normalfont},label={(K\arabic*$^\prime$)}]
\item $x\mapsto \kappa(x,\xi)$ belongs to $BV(\Omega;\R^d)$ for all $\xi\in\R^d$ and there exists a Borel set $N\subset \Omega$ with $\H^{d-1}(N)=0$ such that $\kappa(\cdot,\xi)$ coincides with its lower approximate limit $\kappa^-(\cdot,\xi)$ in $\Omega\setminus N$ for all $\xi\in \R^{d}$\,. \label{ass-k1cc}
\end{enumerate}
Then 
$\phi$ is BV symmetric jointly convex.
\end{prpstn}
\begin{proof}
We may follow the argument of 
Proposition \ref{thm:thm4.13fps} with minor modifications. First, from \ref{ass-k3}, for $\H^{d-1}$-a.e. $x\in\Omega$ and for every $\xi\in\R^d$ a representation as in \eqref{eq:degiorgirepr} holds for $\kappa$, 
where now the functions $a_{j}$ are bounded and, by virtue of \ref{ass-k1cc}, belong to $BV(\Omega;\R^d)$ (see \cite[Remark 2.5]{DCFV2}) 
and 
$
a_{j}={a}^-_{j}.
$
Then, we may define the vector fields $g_{h,j,p}$ as in \ref{eq:w11vectorfiels}, with $b_{j}(x):=a_j(x)+q_jv_j\in \mathcal{A}_x$ for every $j\in\N$ and $\mathcal{A}_x$ as in \eqref{eq:classdomain}, since $b_j\in BV_{\rm loc}(\Omega;\R^d)$. The rest of the argument is exactly the same as for the proof of \eqref{eq:representationformul}, we then omit further details. 
\end{proof}

\subsection{Lower semicontinuity results}

In order to study the lower semicontinuity, firstly we consider the splitting-type model case
\begin{equation}\label{TTT}
\phi(x,r,t,\xi) :=a(x)\kappa(r,t,\xi)\,,
\end{equation}
where $\kappa$ is a symmetric jointly convex function, {according to Definition~\ref{defn:sjconv},} and $a$ is a locally bounded
$BV$ function.
In this case, the proof is obtained by approximating the $BV$ function from below by $W^{1,1}$ functions and using the lower semicontinuity Theorem \ref{thm:gsbdsemicontinuity1} proven for $NA$ symmetric jointly convex integrands.

\begin{prpstn}\label{corcor}
Let $a:\Om\to  [0,+\infty)$ be a locally bounded $BV$ function coinciding with its lower approximate limit $a^-$   and let $\kappa: \R^d\times \R^d\times \R^{d}\to  [0,+\infty)$ be a symmetric jointly convex function. 
Then, for every $(u_n)_n\subset GSBD^p(\Om;\R^d)$ converging in measure to $u\in GSBD^p(\Om;\R^d)$ such that
\begin{equation*}
\sup_{n\in\N}\left[\int_\Om |e(u_n)|^p\,\mathrm{d}x + \H^{d-1}(J_{u_n})\right] <+\infty\,,
\end{equation*}
we have
\begin{equation}\label{eq:e17bbisdd}
\!\!\!\!\!\!\int_{J_u\cap\Om}\!\!\!\!a(x)\kappa(u^+,u^-,\nu_{u})\mathrm{d}\H^{d-1} \leq
\liminf_{n\to +\infty}\!\!\int_{J_{u_n}\cap\Om}\!\!\!a(x)\kappa(u_n^+,u_n^-,\nu_{u_n})\mathrm{d}\H^{d-1}\,.
\end{equation}
\end{prpstn}

\begin{proof}
It suffices to note that by Theorem \ref{t:lecce} the function $a$ is  
lower semicontinuous
with respect to the $1$-capacity. Then by using the approximation result of Lemma \ref{maso1}, we can find an increasing sequence of nonnegative
functions $\{a_h\}\subseteq W^{1,1}(\Omega)$ such that each $a_h$ is approximately continuous $\H^{d-1}$-almost everywhere in $\Omega$, and $a(x)=\sup_{h\in\N}a_h(x)$. Now, we may apply Theorem \ref{thm:gsbdsemicontinuity1} to the sequence $\phi_h(x,r,t,\xi):=a_h(x)\kappa(r,t,\xi)$, whence \eqref{eq:e17bbisdd} follows from Lemma \ref{l:passup}\,. 
\end{proof}

A further step toward the general case is the study of integrands which are lower semicontinuous in $x$ uniformly with respect to the other variables.
For these integrands the following approximation from below holds with functions of the type (\ref{TTT}).

\begin{prpstn}
\label{SCI22}
Let $\phi:\Om\times \R^d\times \R^d\times\R^d\to [0,+\infty)$ be a Borel function such that
\begin{enumerate}[font={\normalfont},label={(LSC)}]
\item given $x_0\in\Om$, for all $\eps>0$
there exists $\delta>0$ such that
\begin{equation*}
\phi(x_{0},r,t,\xi)\leq(1+\eps)\phi(x,r,t,\xi)
\end{equation*}
for all $(x,r,t,\xi)\in\Om\times \R^d\times \R^d\times \R^d$ such that $|x-x_0|<\delta$\,; \label{ass-alpha}
\end{enumerate}
\begin{enumerate}[font={\normalfont},label={(B2$^\prime$)}]
\item for every $x\in\Omega$ the function $\phi(x,\cdot,\cdot,\cdot)$ is symmetric jointly convex\,. \label{ass-beta}
\end{enumerate}
Then for every $j\in\N$ there exist $a_j\in C^\infty_0(\Om;[0,1])$, $a_j(x_j)=1$ for some $x_j\in\Omega$, and $g_j\in  \mathcal{C}(\R^d;\R^d)$ 
such that
\begin{equation}
\phi(x,r,t,\xi)=\sup_{j\in\N}a_j(x)\langle g_j(r)-g_j(t),\xi\rangle^+
\label{eq:supregfunct}
\end{equation}
for all $(x,r,t,\xi)\in \Om \times \R^d\times \R^d\times\R^d$ and $g_j$ satisfies {\rm\ref{ass-C1}}, {\rm\ref{ass-c'}} and {\rm\ref{ass-f}}. 
\end{prpstn}
\begin{proof}
We may follow the argument of the proof of \cite[Proposition~4.5]{DC}, up to replace the assumption of joint convexity by \ref{ass-beta}. The core is the proof of the representation formula
\begin{equation*}
\phi(x,r,t,\xi)= \sup_{G\in\mathcal{G}} G(x,r,t,\xi)\qquad \mbox{ for all }(x,r,t,\xi)\in\Omega\times\R^d\times\R^d\times\R^d\,,
\end{equation*}
where $\mathcal{G}$ is the class of the continuous functions of the form $\varphi(x)\langle h(r)-h(t),\xi\rangle^+$, with $h\in W^{1,\infty}(\R^d;\R^d)$, conservative and $\varphi\in C^\infty_0(\Om;[0,1])$, $\varphi(x)=1$ for some $x\in\Omega$. Then assertion \eqref{eq:supregfunct} follows in a standard way from \cite[Proposition~4.78]{FL}. We omit the details. 
\end{proof}

The lower semicontinuity result for integrands $\phi$ as above is expressed by the following proposition.

\begin{prpstn}\label{t:SCI1}

Let $\phi:\Om\times \R^d\times \R^d\times\R^d\to [0,+\infty)$ be a Borel function such that conditions {\rm\ref{ass-alpha}} and {\rm\ref{ass-beta}} hold. 
Then, for every $(u_n)_n\subset GSBD^p(\Om;\R^d)$ converging in measure to $u\in GSBD^p(\Om;\R^d)$ such that \begin{equation*}
\sup_{n\in\N}\left[\int_\Om |e(u_n)|^p\,\mathrm{d}x + \H^{d-1}(J_{u_n})\right] <+\infty\,,
\end{equation*}
we have
\begin{equation}\label{eq:e17bbis}
\int_{J_u\cap\Om} \phi(x,u^+,u^-,\nu_{u})\,\mathrm{d}\H^{d-1} \leq
\liminf_{n\to +\infty}\int_{J_{u_n}\cap\Om} \phi(x,u_n^+,u_n^-,\nu_{u_n})\,\mathrm{d}\H^{d-1}\,.
\end{equation}
\end{prpstn}

\begin{proof}
By virtue of Proposition~\ref{SCI22}, $\phi$ can be represented according to \eqref{eq:supregfunct} as
\begin{equation*}
\phi(x,r,t,\xi)=\sup_{j\in\N}a_j(x)\langle g_j(r)-g_j(t),\xi\rangle^+\,, 
\end{equation*}
for all $(x,r,t,\xi)\in \Om \times \R^d\times \R^d\times\R^d$, for some sequences of functions $(a_j)_{j\in\N}$ and $(g_j)_{j\in\N}$. Now, since each $a_j$ and $\kappa_j(r,t,\xi):=\langle g_j(r)-g_j(t),\xi\rangle^+$ satisfy the assumptions of Proposition~\ref{corcor}, for every $j\in\N$ we have 
\begin{equation*}
\int_{J_u\cap\Om}a_j(x)\kappa_j(u^+,u^-,\nu_u)\,\mathrm{d}\H^{d-1} \leq
\liminf_{n\to +\infty}\int_{J_{u_n}\cap\Om} a_j(x)\kappa_j(u_n^+,u_n^-,\nu_{u_n})\,\mathrm{d}\H^{d-1}\,.
\end{equation*}
Then, the lower semicontinuity \eqref{eq:e17bbis} is again a consequence of Lemma \ref{l:passup}. 
\end{proof}

In fact, as already remarked in \cite{DC}
for the $GSBV$ setting, assumption {\rm\ref{ass-alpha}} is implied by some conditions which are easier to verify. Thus, the lower semicontinuity result of Proposition~\ref{t:SCI1} still holds if {\rm\ref{ass-alpha}} is replaced by {\rm\ref{ass-b3}} and the following {\rm\ref{ass-alpha1}}, {\rm\ref{ass-alpha2}}.
\begin{prpstn}
\label{SCI2}
Let $\phi:\Om\times \R^d\times \R^d\times\R^d\to [0,+\infty)$ be a Borel function complying with {\rm\ref{ass-beta}}, {\rm\ref{ass-b3}} and such that
\begin{enumerate}[font={\normalfont},label={(C{\arabic*})}]
\item $\phi(\cdot,\cdot,\cdot,\xi)$ is lower semicontinuous on $\Om\times \R^d\times \R^d$
for every $\xi\in\R^d$\,; \label{ass-alpha1}
\item there exists $N\subset\Omega$, with $\H^{d-1}(N)=0$, such that $\phi(x,r,t,\xi)>0$ for all $(x,r,t,\xi)\in(\Om\setminus N)\times \R^d\times \R^d\times(\R^d\setminus \{0\})$\,. \label{ass-alpha2}
\end{enumerate}
Then condition {\rm\ref{ass-alpha}} holds.
\end{prpstn}
\begin{proof}
The argument is essentially that of \cite[Proposition~4.7]{DC}, we then omit the details. 
\end{proof}

The following result still relies on Proposition~\ref{SCI22} (see the argument of \cite[Theorem~4.8]{DC}).

\begin{thrm}\label{mmmdd}
Let $\phi$ be a ${{BV}}$  symmetric jointly convex function satisfying {\rm\ref{ass-alpha2}}. Then the lower semicontinuity {\rm(\ref{eq:e17bbis})} holds.
\end{thrm}

\section{Existence results} \label{sec: existence}

The lower semicontinuity results provided by Theorem~\ref{thm:gsbdsemicontinuity1} and Theorem~\ref{mmmdd}, combined with the compactness theorem in $GSBD$ (recalled with Theorem~\ref{thm:compactness}) allow us to prove the existence of minimizers for functionals of the form 
\begin{equation*}
\mathcal{E}(u):=\int_\Om W(x,e(u))\,\mathrm{d}x + \int_{J_u}\phi(x,u^+,u^-,\nu_u)\,\mathrm{d}\mathcal{H}^{d-1} + \int_\Om \Psi(|u|)\,\mathrm{d}x\,,\quad u\in GSBD^p(\Omega;\R^d)\,,
\end{equation*}
under suitable natural assumptions on the potential $W$ (see $(w1)$--$(w3)$ below).

\begin{thrm}\label{thm:thmexistence}
Let $\Om\subset\R^d$ be open and bounded and let $c>0$. Let $W:\Omega\times\mathbb{M}_{\rm sym}^{d\times d}\to[0,+\infty)$ be such that
\begin{enumerate}
\item[$(w1)$] $W(x,\cdot)$ is convex and lower semicontinuous on $\mathbb{M}_{\rm sym}^{d\times d}$ for a.e. $x\in\Omega$; 
\item[$(w2)$] $W(\cdot, F)$ is measurable on $\Omega$ for every $F\in \mathbb{M}_{\rm sym}^{d\times d}$;
\item[$(w3)$] $W(x,F)\geq c|F|^p$ for all $F\in \mathbb{M}_{\rm sym}^{d\times d}$ for some $p>1$ and for a.e. $x\in\Omega$.
\end{enumerate}
Let $\phi:\Om\times\R^d\times\R^d\times \mathbb{S}^{d-1}\to[c,+\infty)$ be  NA (or BV) symmetric jointly convex and let $\Psi:[0,+\infty)\to[0,+\infty)$ be continuous such that $\displaystyle\lim_{s\to+\infty}\Psi(s)=+\infty$. Then the functional $\mathcal{E}$
admits a minimizer in $GSBD^p(\Om;\R^d)$.
\end{thrm}
\begin{proof}
Let $(u_k)_k\subset GSBD^p(\Om;\R^d)$ be a minimizing sequence. The growth assumption on $W$ and the fact that $\phi\geq c>0$ imply that $u_k$ complies with \eqref{eq:equicoerciv}. Then, by virtue of Theorem~\ref{thm:compactness} we can find $u\in GSBD^p(\Om;\R^d)$ such that $u_k\to u$ a.e. in $\Om$ and $e(u_k)\rightharpoonup e(u)$ weakly in $L^p(\Omega;\mathbb{M}_{\rm sym}^{d\times d})$. From the convexity of $W$ and Fatou's lemma we get
\begin{equation*}
\begin{split}
\int_\Om W(x,e(u))\,\mathrm{d}x & \leq \mathop{\lim\inf}_{k\to+\infty} \int_\Om W(x,e(u_k))\,\mathrm{d}x \\
\int_\Om \Psi(|u|)\,\mathrm{d}x & \leq \mathop{\lim\inf}_{k\to+\infty} \int_\Om \Psi(|u_k|)\,\mathrm{d}x\,.
\end{split}
\end{equation*} 
Finally, the equiboundedness assumption \eqref{eq:equicoerciv} allows us to apply Theorem~\ref{thm:gsbdsemicontinuity1} or Theorem~\ref{mmmdd} to the surface term, thus obtaining
\begin{equation*}
\int_{J_u\cap\Om} \phi(x,u^+,u^-,\nu_{u})\,\mathrm{d}\H^{d-1} \leq
\liminf_{k\to +\infty}\int_{J_{u_k}\cap\Om} \phi(x,u_k^+,u_k^-,\nu_{u_k})\,\mathrm{d}\H^{d-1}\,.
\end{equation*}
Therefore, $u$ is a minimizer of $\mathcal{E}$ and the proof is concluded.
\end{proof}

{\begin{rmrk}
We emphasize the crucial role of the assumption $\phi\geq c>0$ in Theorem~\ref{thm:thmexistence} (combined with $(w3)$), in order to get the compactness in $GSBD^p$ for minimizing sequences $(u_k)$, and then the existence of minimizers in $GSBD^p$. Indeed, such condition allows for a control from below of $\mathcal{E}(u_k)$ with isotropic Griffith-type energies:
\begin{equation*}
c\left(\int_\Om |e(u_k)|^p\,\mathrm{d}x +\mathcal{H}^{d-1}(J_{u_k})\right) + \int_\Om \Psi(|u_k|)\,\mathrm{d}x \leq \mathcal{E}(u_k) <+\infty\,,
\end{equation*}
and the hypotheses of Theorem~\ref{thm:compactness} are satisfied. 
The weaker assumption $\phi\geq0$, in general, would not lead to any control on the surface term $\mathcal{H}^{d-1}(J_{u_k})$, and the compactness in $GSBD^p$ would be not clear.
Nevertheless, in the special case $\phi=\phi(x,r-t,\nu)\geq c |r-t|$, we get compactness in the proper subspace $SBD^p$ (see \cite[Lemma~4.11]{FPS}) and our lower semicontinuity results (Theorem~\ref{thm:gsbdsemicontinuity1} and Theorem~\ref{mmmdd}) can be still applied in some cases (see \cite[Theorem~5.5]{FPS}).
\end{rmrk}}

{
\begin{rmrk}
Several concrete problems in Fracture Mechanics deal with non-convex bulk potentials $W(x,\cdot)$. A reasonable assumption weaker than  $(w1)$ would be the \emph{symmetric quasi-convexity} of the bulk energies (see \cite{Ebobisse}). As far as we know, the only lower semicontinuity result in $GSBD$ under this assumption on the bulk potential is obtained in \cite[Theorem~1.2]{CC22} when restricting to Caccioppoli partitions. The delicate argument therein, based on a blow-up technique, requires a \emph{continuous dependence} of the surface integrand $\phi(x,r-t,\nu)$ on the spatial variable. Nevertheless, if we agree to set the minimization problem in the proper subspace $SBD^p$, 
we can invoke the lower semicontinuity result \cite[Theorem~1.2]{Ebobisse} still for symmetric quasi-convex bulk energies, complying with a stronger growth assumption than $(w3)$. This result could be combined with Theorem~\ref{thm:gsbdsemicontinuity1} or Theorem~\ref{mmmdd} above, up to adding an $L^\infty$-bound on the minimizing sequences $(u_k)$ (see, e.g., \cite[Theorem~3.6]{GZ3} for a similar analysis still in $SBD^p$).
\end{rmrk}}

\section*{Acknowledgements}

The authors are members of Gruppo Nazionale per l’Analisi Matematica, la Probabilità e le loro Applicazioni (GNAMPA) of INdAM. {They acknowledge the INdAM–GNAMPA 2022 Project ``Alcuni problemi associati a funzionali integrali: riscoperta strumenti classici e nuovi sviluppi'', CUPE55F22000270001. The authors are indebted to Francesco Solombrino for many fruitful discussions and for suggesting the key role played by the conservativeness of the approximating vector fields in the proof of Theorem~\ref{chain rule1}. They also gratefully acknowledge the anonymous referees for a careful reading and the useful comments leading to improvement of the manuscript.}

\end{document}